\documentclass[12pt]{amsart}

\usepackage[utf8]{inputenc}
\usepackage{amssymb}
\usepackage{bm}
\usepackage{mathrsfs}
\usepackage{fullpage}
\usepackage[hidelinks]{hyperref}
\usepackage{cleveref}

\newtheorem{thm}{Theorem}[section]
\newtheorem{cor}[thm]{Corollary}
\newtheorem{lem}[thm]{Lemma}
\newtheorem{prop}[thm]{Proposition}
\theoremstyle{definition}
\newtheorem{defn}[thm]{Definition}
\theoremstyle{remark}

\crefname{thm}{Theorem}{Theorems}
\crefname{cor}{Corollary}{Corollaries}
\crefname{lem}{Lemma}{Lemmas}
\crefname{prop}{Proposition}{Propositions}
\crefname{defn}{Definition}{Definitions}
\crefname{rem}{Remark}{Remarks}
\numberwithin{equation}{section}

\newcommand{\norm}[1]{\left\Vert#1\right\Vert}
\newcommand{\abs}[1]{\left\vert#1\right\vert}
\newcommand{\bF}{\mathbb F}
\newcommand{\bQ}{\mathbb Q}
\newcommand{\bZ}{\mathbb Z}
\newcommand{\bN}{\mathbb N}
\newcommand{\bR}{\mathbb R}

\newcommand{\cO}{\mathcal O}
\newcommand{\cL}{\mathcal L}
\newcommand{\cN}{\mathcal N}

\newcommand{\ord}{\operatorname{ord}_p}

\newcommand{\CT}{\operatorname{CT}}
\newcommand{\SC}{\operatorname{SC}}

\begin{document}

\title{On $p$-adic asymptotic formulas for slope factors of $L$-functions of toric exponential sums}
\author{Chunlin Wang}
\address{School of Mathematical Sciences, Sichuan Normal University, Chengdu 610066, China}
\email{c-l.wang@outlook.com}
\subjclass[2020]{11S40, 11T23, 14G10}
\keywords{exponential sums, nuclear $\sigma$-modules, Hodge--Newton decomposition, Schur complements, Dwork trace formula}

\begin{abstract}
Adolphson and Sperber expressed the unique unit root of the $L$-function of
a toric exponential sum in terms of $A$-hypergeometric functions, whereas
reciprocal zeros and poles of positive slope remain much less understood.
We study complete rational factors of arbitrary prescribed slope in families
satisfying a chain-level ordinary condition.  Our main tool is a new
contraction map arising from a twisted block LU factorization of the
Frobenius matrix of a nuclear $\sigma$-module.  Its fixed point produces an
upper triangularization and the associated slope filtration, while explicit
Schur-complement quotients of finite Frobenius iterates approximate the
diagonal slope blocks with geometrically decreasing error.  At a periodic
fibre, successive cancellation in an ordered product of these approximations
yields a quotient formed from two finite Frobenius iterates, and the characteristic polynomial
of this quotient converges to the complete slope polynomial.  Combined with the relative
Dwork trace formula, this yields uniform $p$-adic asymptotic formulas for
the complete slope factors of the $L$-functions of toric exponential sums.  When the
relevant slope blocks have rank one, the Schur complements reduce to ratios
of leading principal minors and give asymptotic formulas for the individual
reciprocal zeros and poles after cancellation.
\end{abstract}

\maketitle

\section{Introduction}\label{sec:introduction}

Let $p$ be an odd prime, let $q$ be a power of $p$, and fix a primitive
$p$-th root of unity $\zeta_p\in\overline{\bQ}_p$.  Write
\[
  \Lambda=(\Lambda_1,\ldots,\Lambda_m),
  \qquad X=(X_1,\ldots,X_n),
\]
and consider a Laurent polynomial
\begin{equation}\label{eq:laurent-polynomial}
  f(\Lambda,X)
  =\sum_{u\in\bN^m,\,v\in\bZ^n}
    a_{u,v}\Lambda^uX^v
  \in\bF_q[\Lambda][X_1^{\pm1},\ldots,X_n^{\pm1}].
\end{equation}
For $\lambda\in\overline{\bF}_q^m$, let $\bF_q(\lambda)$ be the field
generated by its coordinates, put
\[
  d(\lambda)=[\bF_q(\lambda):\bF_p],
\]
and denote by $\bF_q(\lambda;k)$ its extension of degree $k$.  With
$T_{\lambda,k}$ denoting the trace to $\bF_p$, set
\[
  S_{\lambda,k}
  =\sum_{x\in(\bF_q(\lambda;k)^*)^n}
    \zeta_p^{T_{\lambda,k}(f(\lambda,x))}
\]
and
\begin{equation}\label{eq:L-function-intro}
  L_\lambda(T)
  =\exp\left(\sum_{k\geq1}\frac{S_{\lambda,k}}{k}T^k\right).
\end{equation}
This $L$-function is rational, and its reciprocal zeros and poles are
algebraic integers.  Their complex and $p$-adic properties have been studied
extensively; see \cite{De80,AS89,DL91,W93,W04}.

Beyond their valuations, one would like to understand how these reciprocal
zeros and poles vary with $\lambda$.  Dwork expressed the unit root of the
Legendre family in terms of a $p$-adic hypergeometric function \cite{Dw69}
and obtained an analogous formula for the Kloosterman family using a
$p$-adic Bessel function \cite{Dw74}.  Sperber extended the latter work to
higher-dimensional Kloosterman sums \cite{Sp80}, and Adolphson and Sperber
gave an $A$-hypergeometric formula for the unique unit root of a broad class
of toric exponential sums \cite{AS12}.  
Related work on unit-root
$L$-functions appears in \cite{HS14,YZ26}; see
\cite{Dw73,W99,W00a,W00b} for the surrounding theory.

Formulas for reciprocal zeros and poles of positive slope are much less developed,
apart from special cases such as the Airy family studied by Zhang \cite{Zh}, and the hyperkloosterman family stuidied by Adolphson, Sperber and Zhang \cite{ASZ}.
The aim of this paper is to construct explicit, uniform $p$-adic approximations to complete rational factors of prescribed slope in toric
families that are ordinary on the chain level through those slopes. The key ingredient is a contraction map obtained from a twisted block LU factorization of Frobenius. 

Fixed-point methods also appear in Wan's work
\cite{W93,W99,W00a}, where they are used to construct slope filtrations,
although the underlying contraction is not isolated in the explicit form
needed here. In the rank-one unit-root setting, Adolphson and Sperber
\cite{AS12} used an explicit contraction in their construction of the
unique unit root. Our construction combines these two perspectives,
extending the contraction to block changes of basis so as to treat entire
slope blocks. In particular, we recover the formula of Adolphson and Sperber for the unique unit root of a family (See \Cref{subsec:as-unit-root-comparison}).

The argument has two parts.  First, a twisted block LU factorization
produces a contraction on a space of formal basis changes for a nuclear
$\sigma$-module.  Its fixed point gives the slope filtration.  More
importantly for the present purpose, each finite approximation to a graded
Frobenius block is an explicit quotient of Schur complements of two
successive iterated Frobenius matrices.  At a periodic fibre, the intermediate
Schur factors in an ordered product at consecutive approximation levels
cancel successively, leaving a quotient involving only two finite Frobenius
iterates.  Second, the relative Dwork trace
formula translates the resulting characteristic-polynomial approximations
into asymptotic formulas for the slope factors of the family \eqref{eq:L-function-intro}. 

The contraction theorem applies to both finite-rank and countable-rank
nuclear $\sigma$-modules and therefore suggests applications at both the
cohomological and chain levels.  In this paper, the applications in
\Cref{sec:dwork,sec:application} are carried out at the chain level.  We do
not pursue a general study of cohomology-level applications, which would
lead beyond the main scope of the paper.  Instead,
\Cref{sec:differential-interpretation} presents a simple one-parameter
finite-rank differential-module illustration, motivated by possible
cohomological applications, in which additional structure
identifies the limiting graded Frobenius blocks directly from the upper
factor in the block LU decomposition of a constant-gauge transform of a
normalized horizontal solution matrix.  This gives exact
graded Frobenius formulas and, in rank one, an exact root formula,
illustrating the potential of the contraction construction at the
cohomological level.

We now give a concise statement of the abstract construction.  The complete
definitions and proofs are in \Cref{sec:sigma-modules}.  Let $(M,\phi)$ be a
nuclear $\sigma$-module over a complete normed coefficient ring $A_0$, and
let $\xi\in A_0$ satisfy $0<\rho:=\abs{\xi}<1$.  Suppose that, in a formal
basis grouped into finite blocks of sizes $h_0,h_1,\ldots$, the Frobenius
matrix has the form
\[
  B=(\xi^jB_{ij})_{i,j\geq0}.
\]
Put
\[
  H_j=\sum_{r=0}^jh_r,
  \qquad
  N_j=\sum_{r=0}^jrh_r,
\]
and let $\Delta_j(B)$ be the leading principal minor of order $H_j$.
Set $\Delta_{-1}(B)=1$.  The global matrix hypothesis used below is
\begin{equation}\label{eq:intro-ordinary-minors}
  \Delta_j(B)=\xi^{N_j}u_j,
  \qquad u_j\in A_0^\times,
  \qquad 0\leq j\leq k.
\end{equation}
It implies fibrewise ordinarity through slope $k$.  For the coefficient
ring used in the Dwork construction below, the Teichm\"uller specializations
detect units, so the two conditions are equivalent.

Let $\mathfrak I_k$ be the space of lower block unitriangular matrices whose
off-diagonal blocks are supported in block columns $0,\ldots,k$.  For
$I\in\mathfrak I_k$, normalized block elimination gives a unique
factorization
\[
  BI^\sigma=J(I)K(I),
\]
where $J(I)\in\mathfrak I_k$ and the first $k+1$ block columns of $K(I)$
are upper block triangular.

For a block matrix $A$ whose leading principal block through index $i-1$
is invertible over $A_0[\xi^{-1}]$, define its $i$th Schur complement by
\[
  \SC_i(A)
  =A_{ii}-A_{i,<i}A_{<i,<i}^{-1}A_{<i,i},
\]
where $A_{<i,<i}$ is the leading $H_{i-1}\times H_{i-1}$ block; set
$\SC_0(A)=A_{00}$.  Thus $\SC_i(A)$ depends only on the leading
$H_i\times H_i$ truncation of $A$.
Write $\Delta_i(A)$ for the determinant of this leading truncation, set
$\Delta_{-1}(A)=1$, and put
\[
  \widehat\Delta_i(A)
  :=\det\SC_i(A)
  =\frac{\Delta_i(A)}{\Delta_{i-1}(A)}.
\]

\begin{thm}[Block contraction and graded-block approximation]
\label{thm:intro-main}
Under \eqref{eq:intro-ordinary-minors}, the map $\alpha(I)=J(I)$ is a
contraction with constant $\rho$ and hence has a unique fixed point $P$.
If $Q=K(P)$, then
\[
  BP^\sigma=PQ,
\]
the first $k+1$ block columns of $Q$ are upper block triangular, and their
diagonal blocks have the form
\[
  D_i=Q_{ii}=\xi^iC_{ii},
  \qquad C_{ii}\in\operatorname{GL}_{h_i}(A_0),
  \qquad 0\leq i\leq k.
\]
Starting with $I_0=1$, define,
for $a\geq1$,
\[
  I_a=\alpha(I_{a-1}),
  \qquad BI_{a-1}^\sigma=I_aK_a,
  \qquad D_i^{(a)}=(K_a)_{ii}.
\]
Put
\[
  B_0=1,
  \qquad
  B_a=BB^\sigma\cdots B^{\sigma^{a-1}}
  \quad(a\geq1).
\]
Then $\xi^{-ai}\SC_i(B_a)\in\operatorname{GL}_{h_i}(A_0)$ and
\begin{equation}\label{eq:intro-schur-quotient}
  D_i^{(a)}
  =\SC_i(B_a)\bigl(\SC_i(B_{a-1})^\sigma\bigr)^{-1}.
\end{equation}
Moreover,
\[
  \xi^{-i}(D_i^{(a)}-D_i)\in\xi^aM_{h_i}(A_0),
  \qquad
  \norm{\xi^{-i}(D_i^{(a)}-D_i)}\leq\rho^a.
\]
If $h_i=1$, then
\[
  D_i^{(a)}
  =\frac{\widehat\Delta_i(B_a)}
         {\widehat\Delta_i(B_{a-1})^\sigma},
  \qquad
  \abs{\xi^{-i}\left(
   \frac{\widehat\Delta_i(B_a)}
        {\widehat\Delta_i(B_{a-1})^\sigma}-D_i\right)}
  \leq\rho^a.
\]
\end{thm}

The fixed-point identity $BP^\sigma=PQ$ also recovers the ordinary
Hodge--Newton filtration: the successive spans of the first $i+1$ blocks
of $\bm eP$ form an increasing $\phi$-stable filtration, transversal to
the original basis filtration, and $D_i$ is the Frobenius block on its
$i$th graded piece.  If $x$ is a $\sigma$-periodic specialization of
degree $d$, put
\[
  \Phi_{i,x}
  =D_i(x)D_i^\sigma(x)\cdots D_i^{\sigma^{d-1}}(x),
  \qquad
  F_{i,x}(T)=\det(1-T\Phi_{i,x}).
\]
Then $F_{i,x}(T)$ is the complete slope-$i$ factor of the fibrewise
characteristic series.

\begin{thm}[Fibrewise characteristic-polynomial approximation]
\label{thm:intro-fibre-approximation}
Let $x$ be a $\sigma$-periodic specialization of degree $d$, and put
\begin{equation}\label{eq:intro-periodic-schur-quotient}
 \begin{aligned}
  \Phi_{i,x}^{(a)}
  &:=D_i^{(a+d-1)}(x)
   (D_i^{(a+d-2)})^\sigma(x)\cdots
   (D_i^{(a)})^{\sigma^{d-1}}(x)\\
  &=\SC_i(B_{a+d-1})(x)\SC_i(B_{a-1})(x)^{-1}.
 \end{aligned}
\end{equation}
Then
\[
  \norm{\xi^{-id}(\Phi_{i,x}^{(a)}-\Phi_{i,x})}
  \leq\rho^a.
\]
For $G(T)=\sum_{m=0}^{h_i}g_mT^m\in K_x[T]$, define the weighted
Gauss norm
\[
  \norm{G}_{i,d}
  :=\max_{0\leq m\leq h_i}\rho^{-idm}\abs{g_m}
  =\norm{G(\xi^{-id}T)}_{\mathrm G}.
\]
Then
\begin{equation}\label{eq:intro-periodic-polynomial-error}
  \norm{\det(1-T\Phi_{i,x}^{(a)})-F_{i,x}(T)}_{i,d}
  \leq\rho^a.
\end{equation}
If $h_i=1$, writing $\gamma_{i,x}=\Phi_{i,x}$ and
\[
  \gamma_{i,x}^{(a)}
  :=\frac{\widehat\Delta_i(B_{a+d-1})(x)}
          {\widehat\Delta_i(B_{a-1})(x)},
\]
one has
\[
  \abs{\xi^{-id}
    (\gamma_{i,x}^{(a)}-\gamma_{i,x})}
  \leq\rho^a.
\]
\end{thm}

The new content is concentrated in \Cref{thm:intro-main} and
\Cref{thm:intro-fibre-approximation}: the former gives the finite
Schur-complement formula \eqref{eq:intro-schur-quotient} for the graded
blocks, while the latter turns the ordered periodic product into the
two-iterate quotient \eqref{eq:intro-periodic-schur-quotient} and proves
the characteristic-polynomial estimate
\eqref{eq:intro-periodic-polynomial-error}.  In higher rank the full
matrix-valued Schur complement is needed to recover all characteristic
coefficients; leading principal minors determine only the determinant and
hence the top coefficient.  In rank one these two levels coincide and give
an approximation to the individual fibrewise eigenvalue.

We apply the theorem to the dual relative Dwork operator of
\eqref{eq:laurent-polynomial}.  Let $\Delta$ be the convex hull of the
origin and the $X$-support of $f$, let $w$ be its weight function, and choose
the least integer $D\geq1$ such that
$w(C(\Delta)\cap\bZ^n)\subset D^{-1}\bZ_{\geq0}$.  The relative trace
formula expresses the fibrewise $L$-function as
\[
  L_\lambda(T)^{(-1)^{n-1}}
  =\prod_{\ell=0}^{n}
    D_d(\lambda,p^{d\ell}T)^{
      (-1)^\ell\binom{n}{\ell}},
  \qquad d=d(\lambda),
\]
where $D_d(\lambda,T)$ is the Fredholm determinant of the local dual Dwork
operator at the chain level.

Suppose that the dual Dwork module is ordinary through slope $k/D$, and let
$F_{u,\lambda}(T)$ denote the characteristic polynomial of its complete
slope-$u/D$ Frobenius block.  Then, for $0\leq r\leq k$,
\[
  [L_\lambda(T)^{(-1)^{n-1}}]_{r/D}
  =\prod_{\ell=0}^{\min(n,\lfloor r/D\rfloor)}
    F_{r-\ell D,\lambda}(p^{d\ell}T)^{
      (-1)^\ell\binom{n}{\ell}},
\]
where $[L_\lambda(T)^{(-1)^{n-1}}]_{r/D}$ denote the slope-$r/D$ part of $L_\lambda(T)$. 
For $0\leq r\leq k$ with $r<D$, only the term $\ell=0$ occurs, and hence
\[
  [L_\lambda(T)^{(-1)^{n-1}}]_{r/D}=F_{r,\lambda}(T).
\]
Thus the ordinary slope factor below one is controlled by a single
chain-level Frobenius block.  The new point is that each
$F_{u,\lambda}(T)$ has an approximation expressed entirely in terms of
finite Dwork iterates.  Let $B$ be the dual
relative Dwork matrix.  For $a\geq1$ and $d=d(\lambda)$, define
\[
  F_{u,\lambda}^{(a)}(T)
  =\det\!\left(
    1-T\SC_u(B_{a+d-1})(\lambda)
       \SC_u(B_{a-1})(\lambda)^{-1}
  \right).
\]
By \eqref{eq:intro-periodic-polynomial-error},
\[
  \norm{F_{u,\lambda}^{(a)}-F_{u,\lambda}}_{u,d}
  \leq\abs{\xi}^a,
\]
uniformly in $\lambda$.  Consequently,
\[
  \mathcal R_{r,\lambda}^{(a)}(T)
  =\prod_{\ell=0}^{\min(n,\lfloor r/D\rfloor)}
   F_{r-\ell D,\lambda}^{(a)}
   (p^{d\ell}T)^{(-1)^\ell\binom n\ell}
\]
has the following precise interpretation: when the terms with positive and
negative exponents are placed in the numerator and denominator,
respectively, the two resulting products converge coefficientwise to the
corresponding numerator and denominator products in the trace-formula
expression for the complete slope-$r/D$ factor of
$L_\lambda(T)^{(-1)^{n-1}}$, before any cancellation.  In particular, an
ordinary slope below one
is approximated by the characteristic polynomial of a single explicit
Schur-complement quotient.

When the slope-$u/D$ block has rank one, the Schur complements are scalars and
\[
  \SC_u(B_a)=\widehat\Delta_u(B_a).
\]
Thus the corresponding reciprocal Frobenius root is
\[
  \gamma_{u,\lambda}^{(a)}
  :=\frac{\widehat\Delta_u(B_{a+d-1})(\lambda)}
          {\widehat\Delta_u(B_{a-1})(\lambda)},
  \qquad
  \abs{\xi^{-ud}
    (\gamma_{u,\lambda}^{(a)}-\gamma_{u,\lambda})}
  \leq\abs{\xi}^a.
\]
After the trace-formula factors of opposite signs have been cancelled,
these rank-one formulas describe the surviving reciprocal zeros and poles.
In higher rank, the construction instead determines the complete slope
polynomial, without selecting a canonical ordering of its individual
reciprocal roots.  No smoothness or nondegeneracy hypothesis enters the
chain-level argument beyond the standing Dwork assumptions specified in
\Cref{sec:dwork}.

The paper is organized as follows.  In \Cref{sec:sigma-modules} we develop
the contraction and its fixed-point triangularization, identify the
fibrewise slope factors from the diagonal blocks, and then derive their
Schur-complement and characteristic-polynomial approximations.  In
\Cref{sec:dwork} we recall the relative Dwork
operator, the trace formula, and its dual form.  In
\Cref{sec:application} we apply the Schur formulas to obtain explicit
finite-iterate approximations first to the Dwork graded blocks, then to
their fibrewise slope polynomials and the resulting complete rational
slope factors of the $L$-functions; leading principal minors give
individual asymptotic root formulas in rank one.
Finally, \Cref{sec:differential-interpretation} gives a one-parameter
finite-rank differential-module illustration motivated by possible
cohomological applications: a
constant-gauge transform of a normalized horizontal solution matrix
realizes the contraction fixed point through block LU factorization, while
its upper factor gives exact formulas for the graded Frobenius blocks and
their periodic characteristic polynomials.

\section{The contraction construction for ordinary nuclear
\texorpdfstring{$\sigma$}{sigma}-modules}
\label{sec:sigma-modules}

This section gives the precise abstract setting and proves the contraction
and approximation results stated in the Introduction.

\subsection{Nuclear \texorpdfstring{$\sigma$}{sigma}-modules and ordinarity}

Let $K$ be a complete non-Archimedean valued field with valuation ring
$\mathcal O_K$.  Let $A_0$ be a commutative $\mathcal O_K$-algebra endowed
with a norm $\abs{\cdot}$ extending the norm on $\mathcal O_K$.  Here a
norm means a function $A_0\to\bR_{\geq0}$ such that
\[
 \begin{gathered}
  \abs{a}=0\Longleftrightarrow a=0,
  \qquad \abs{1}=1,\\
  \abs{a+b}\leq\max\{\abs{a},\abs{b}\},
  \qquad \abs{ab}\leq\abs{a}\abs{b}
 \end{gathered}
\]
for all $a,b\in A_0$.  We assume that $\abs{a}\leq1$ for every
$a\in A_0$ and that $A_0$ is complete for this norm.  Fix
$\xi\in\mathcal O_K$ whose image in $A_0$ is a non-zero-divisor and such
that
\[
  0<\rho:=\abs{\xi}<1.
\]
We further assume that the norm topology of $A_0$ agrees with its
$\xi$-adic topology.  Let $\sigma:A_0\to A_0$ be an isometric endomorphism
whose restriction to $\mathcal O_K$ is an automorphism, and suppose that
$\sigma(\xi)=\xi$.

A Banach $A_0$-module with a finite or countable \emph{formal basis}
$\bm e=(e_s)_{s\in S}$, where $S=\{1,\ldots,N\}$ or $S=\bN_{\geq1}$,
is a module of the form
\begin{equation}\label{eq:formal-basis-module}
  M=\left\{\sum_{s\in S}a_se_s:a_s\in A_0\right\},
  \qquad
  \norm{\sum_{s\in S}a_se_s}=\sup_{s\in S}\abs{a_s}.
\end{equation}
Thus the coefficients are bounded, because every element of $A_0$ has norm
at most one, but they need not tend to zero.  This distinguishes a formal
basis from an orthonormal basis.

\begin{defn}\label{def:nuclear-sigma-module}
A \emph{nuclear $\sigma$-module} over $A_0$ is a pair $(M,\phi)$, where
$M$ has a finite or countable formal basis and $\phi:M\to M$ is a
$\sigma$-linear map.  In the countably infinite case we require
\begin{equation}\label{eq:nuclear-condition}
  \lim_{s\to\infty}\norm{\phi(e_s)}=0.
\end{equation}
In finite rank the nuclearity condition is automatic.
The map $\phi$ is called Frobenius.
\end{defn}

Because the norm and $\xi$-adic topologies agree, in countably infinite
rank \eqref{eq:nuclear-condition} is equivalent to requiring that, for
every $r\geq1$, all but finitely many $\phi(e_s)$ belong to $\xi^rM$.
The latter condition is automatic in finite rank.  In either case, let
$d_r$ be the least nonnegative integer such that
\[
  \phi(e_s)\in\xi^rM\qquad(s\in S,\ s>d_r),
\]
and put $d_0=0$ and $h_r=d_{r+1}-d_r$.  The sequence
$h(\bm e)=(h_0,h_1,\ldots)$ is the \emph{basis sequence}.  If $M(r)$ is
the formal submodule with basis $e_s$ for $s\in S$ and $s>d_r$, then
\[
 M=M(0)\supset M(1)\supset\cdots,
 \qquad
 \phi(M(r))\subseteq\xi^rM,
 \qquad
 \operatorname{rank}_{A_0}M(r)/M(r+1)=h_r.
\]
We call this the \emph{basis filtration}.

Order the basis by this filtration and group it into blocks $\bm e_r$ of
size $h_r$.  If $\phi(\bm e)=\bm eB$, then
\begin{equation}\label{eq:abstract-block-matrix}
  B=(\xi^jB_{ij})_{i,j\geq0},
  \qquad
  B_{ij}\in M_{h_i\times h_j}(A_0).
\end{equation}
Put
\[
  H_j=\sum_{r=0}^jh_r,
  \qquad
  N_j=\sum_{r=0}^jrh_r.
\]
The \emph{basis polygon} $P(\bm e)$ is the lower convex polygon joining
$(0,0)$ and the points $(H_j,N_j)$.  Its side of index $j$ has slope $j$
and horizontal length $h_j$.
The terms basis sequence, basis filtration, and basis polygon follow
\cite[Definitions~4.1 and~4.2]{W00a}.

To describe fibrewise ordinarity, fix a collection $\mathcal X$ of bounded
$\mathcal O_K$-algebra homomorphisms
\begin{equation}\label{eq:specializations}
  x:A_0\longrightarrow\mathcal O_{K_x},
\end{equation}
where $K_x/K$ is a finite valued extension.  All these fields have the same
characteristic and residue characteristic as $K$, and $x(\xi)=\xi$.  The
integral fibre relative to $\bm e$ is the formal base change
\[
  M_x^\circ
  =\left\{
    \sum_{s\in S}b_se_{s,x}:b_s\in\mathcal O_{K_x}
   \right\},
  \qquad
  \norm{\sum_{s\in S}b_se_{s,x}}=\sup_{s\in S}\abs{b_s},
\]
with coefficientwise specialization $M\to M_x^\circ$.  Its associated
$K_x$-Banach space is
\[
  M_x=M_x^\circ\otimes_{\mathcal O_{K_x}}K_x
  =\left\{
    \sum_{s\in S}b_se_{s,x}:b_s\in K_x,\
      \sup_{s\in S}\abs{b_s}<\infty
   \right\}.
\]
The construction of $M_x^\circ$ is relative to the chosen formal basis.  If,
in addition, $K$ is discretely valued and $A_0$ is Noetherian, the natural
specialization map induces an isomorphism
\[
  M\otimes_{A_0,x}\mathcal O_{K_x}\xrightarrow{\;\sim\;}M_x^\circ,
\]
so in this case the formal fibre agrees with the algebraic base change.  In
general, however, this map need not be an isomorphism.  We therefore use the
basis-relative formal fibre defined above; none of the arguments below
requires an identification with the algebraic tensor product.

The specialization $x$ is \emph{$\sigma$-periodic of degree $d$} if $d$ is
the least positive integer such that
\begin{equation}\label{eq:periodic-specialization}
  x\circ\sigma^d=x.
\end{equation}
Then $\phi^d$ induces a $K_x$-linear endomorphism of $M_x$ whose matrix is
\[
  B_{[x]}=B(x)B^\sigma(x)\cdots B^{\sigma^{d-1}}(x).
\]
Normalize the additive valuation $v_x$ by $v_x(\xi)=1$.  The Newton polygon
at $x$ is the Newton polygon of $\det(1-TB_{[x]})$, with $v_x$ divided by
$d$.  Following \cite[Definition~4.8]{W00a}, the basis $\bm e$ is
\emph{ordinary up to the slope $k$ side} relative to $\mathcal X$ if this
Newton polygon agrees with $P(\bm e)$ through all sides of slope $\le k$, 
for every $x\in\mathcal X$.

For a block matrix $C$, let $\Delta_j(C)$ denote its leading principal
minor of order $H_j$, and set $\Delta_{-1}(C)=1$.

\begin{defn}[Leading-minor hypothesis]\label{def:block-ordinary}
The Frobenius matrix $B$ satisfies the \emph{leading-minor hypothesis
through slope $k$} if
\begin{equation}\label{eq:ordinary-leading-minors}
  \Delta_j(B)=\xi^{N_j}u_j,
  \qquad u_j\in A_0^\times,
  \qquad 0\leq j\leq k.
\end{equation}
\end{defn}

The literature formulates ordinarity in terms of equality of the basis and
Newton polygons, as above.  The leading-minor hypothesis is the global
coefficient-ring criterion needed for the contraction construction; it
implies that the chosen basis is ordinary up to the slope-$k$ side at every
periodic fibre.  The converse holds when the chosen specializations detect
units:
\[
  \abs{x(g)}=1\text{ for every }x\in\mathcal X
  \quad\Longrightarrow\quad g\in A_0^\times.
\]
This property holds for the Teichm\"uller specializations of the ring
$\cO^\circ$ of power-bounded elements in the Tate algebra introduced in
\Cref{sec:dwork}; see \Cref{prop:ordinarity-criterion}.

The abstract contraction argument does not require characteristic zero.
The common base field in \eqref{eq:specializations} is used only to make the
valuations and residue characteristics of the fibres compatible.  The
Dwork application is, of course, over a characteristic-zero $p$-adic
field.

For the remainder of this section, fix $k\geq0$ and assume
\eqref{eq:ordinary-leading-minors}.

\subsection{The block contraction}

Let $\mathfrak I_k$ be the set of lower block
unitriangular matrices $I=(I_{ij})_{i,j\geq0}$ such that
\[
 I_{ii}=1,
 \qquad I_{ij}=0\ (i<j),
 \qquad I_{ij}=0\ (j>k,\ i\ne j),
\]
and all entries lie in $A_0$.  Equip it with
\begin{equation}\label{eq:I-space-norm}
  \norm{I-I'}=\sup_{\substack{i>j\\0\leq j\leq k}}
     \norm{I_{ij}-I'_{ij}}.
\end{equation}
The space is complete because it is a closed product of copies of $A_0$.
Moreover, every row of $I$ has finite support away from the diagonal.  The
same is true of $I^{-1}$.  Indeed, $I-1$ is nilpotent of exponent at most
$k+2$, because every nonzero off-diagonal column has index at most $k$; thus
$I^{-1}=1-(I-1)+\cdots+(-1)^{k+1}(I-1)^{k+1}$.  Hence
multiplication by $I$ and $I^{-1}$ is defined on every bounded coefficient
sequence, and every $I\in\mathfrak I_k$ gives a legitimate formal change of
basis.

\begin{lem}\label{lem:minors-preserved}
Assume \eqref{eq:ordinary-leading-minors}.  For every
$I\in\mathfrak I_k$ and $0\leq j\leq k$,
\[
  \Delta_j(BI^\sigma)=\xi^{N_j}u_j(I),
  \qquad u_j(I)\in A_0^\times.
\]
\end{lem}

\begin{proof}
Recall that
\[
  B=(\xi^jB_{ij})_{i,j\geq0}.
\]
Thus, writing \(B_{\ast r}=(B_{ir})_{i\geq0}\), the actual \(r\)th block
column of $B$ is $\xi^rB_{\ast r}$.  Since $I$ is lower block
unitriangular, for $0\leq r\leq j$ the $r$th block column of $BI^\sigma$ is
\begin{equation}\label{eq:column-operation}
\begin{aligned}
  (BI^\sigma)_{\ast r}
   &=\xi^rB_{\ast r}
     +\sum_{\ell>r}\xi^\ell
        B_{\ast\ell}I_{\ell r}^\sigma  \\
   &=\xi^r\left(
        B_{\ast r}
        +\sum_{\ell>r}\xi^{\ell-r}
           B_{\ast\ell}I_{\ell r}^\sigma
      \right).
\end{aligned}
\end{equation}
After extracting \(\xi^r\) from each of the \(h_r\) columns in the
\(r\)th block column, every term added to \(B_{\ast r}\) is divisible
by \(\xi\), since \(\ell-r\geq1\).

Restricting the first \(j+1\) block columns in
\eqref{eq:column-operation} to the first \(H_j\) rows and taking
determinants therefore gives
\[
  \xi^{-N_j}\Delta_j(BI^\sigma)
  \equiv
  \xi^{-N_j}\Delta_j(B)
  \pmod{\xi A_0}.
\]
Indeed, extracting \(\xi^r\) from the \(h_r\) scalar columns belonging
to block column \(r\), for \(0\leq r\leq j\), produces the total factor
\[
  \prod_{r=0}^j\xi^{rh_r}=\xi^{N_j}.
\]

By the leading-minor hypothesis,
\[
  \xi^{-N_j}\Delta_j(B)=u_j\in A_0^\times.
\]
Consequently,
\[
  \xi^{-N_j}\Delta_j(BI^\sigma)
   =u_j+\xi a
   =u_j(1+\xi a')
\]
for some \(a,a'\in A_0\).  Since \(\lvert\xi\rvert<1\) and every
element of \(A_0\) has norm at most one, \(1+\xi a'\) is a unit, with
inverse given by its convergent geometric series.  Hence
\[
  \Delta_j(BI^\sigma)=\xi^{N_j}u_j(I)
\]
for some \(u_j(I)\in A_0^\times\), as required.
\end{proof}

The next lemma is the \emph{normalized block elimination} used in
\Cref{thm:intro-main}.  The normalization requires the lower factor to have
identity diagonal blocks and extracts $\xi^j$ from the $j$th block column of
the upper factor.

\begin{lem}[Normalized block elimination]\label{lem:block-lu}
Assume \eqref{eq:ordinary-leading-minors}.  For every
$I\in\mathfrak I_k$ there are unique matrices
\[
  J=J(I)\in\mathfrak I_k,
  \qquad K=K(I),
\]
such that
\[
  BI^\sigma=JK,
\]
and the first $k+1$ block columns of $K$ are block upper triangular.  More
precisely, for $0\leq j\leq k$,
\[
  K_{ij}=0\quad(i>j),
  \qquad
  K_{ij}=\xi^jC_{ij}\quad(i\leq j),
\]
where $C_{ij}$ has entries in $A_0$ and
$C_{jj}\in\operatorname{GL}_{h_j}(A_0)$.
\end{lem}

\begin{proof}
We perform block Gaussian elimination from left to right by induction.  For
the initial step $h=0$, set, for $0\leq j\leq k$,
\[
  C_{0j}
  =B_{0j}
   +\sum_{\ell>j}\xi^{\ell-j}
      B_{0\ell}I_{\ell j}^{\sigma}.
\]
The sum converges because its $\ell$th term has norm at most
$\rho^{\ell-j}$.  Then $K_{0j}=\xi^jC_{0j}$ is the $(0,j)$ block of
$BI^\sigma$.  Moreover,
\[
  \det C_{00}
  =\Delta_0(BI^\sigma)\in A_0^\times
\]
by \Cref{lem:minors-preserved}, since $N_0=0$.  Hence $C_{00}$ is
invertible.  Taking $J_{00}=1$, define, for $i>0$,
\[
  J_{i0}
  =\left(
      B_{i0}
      +\sum_{\ell>0}\xi^\ell
         B_{i\ell}I_{\ell0}^{\sigma}
    \right)C_{00}^{-1}.
\]
This gives
\[
  (BI^\sigma)_{i0}=J_{i0}C_{00}=(JK)_{i0},
\]
and therefore eliminates the blocks below the first pivot.  This establishes
the initial step of the induction.

Now let \(1\leq h\leq k\), and suppose that block columns
\(0,\ldots,h-1\) have already been eliminated.  For \(h\leq j\leq k\),
comparison of the \((h,j)\) block of \(BI^\sigma=JK\), after dividing the
\(j\)th block column by \(\xi^j\), gives
\begin{equation}\label{eq:C-recursion}
  C_{hj}
  =B_{hj}
   +\sum_{\ell>j}\xi^{\ell-j}
      B_{h\ell}I_{\ell j}^{\sigma}
   -\sum_{\ell<h}J_{h\ell}C_{\ell j}.
\end{equation}
The first sum converges because its \(\ell\)th term has norm at most
\(\rho^{\ell-j}\), while the second sum is finite.  Thus every entry of
\(C_{hj}\) belongs to \(A_0\).

We next show that the pivot \(C_{hh}\) is invertible.  The leading
\(H_h\times H_h\) block of \(J\) is lower block unitriangular, whereas the
corresponding block of \(K\) is upper block triangular with diagonal blocks
\(\xi^rC_{rr}\), \(0\leq r\leq h\).  Consequently,
\[
  \xi^{-N_h}\Delta_h(BI^\sigma)
  =\prod_{r=0}^h\det C_{rr}.
\]
The left-hand side is a unit by \Cref{lem:minors-preserved}.  Since the
previous pivots are invertible by induction, \(\det C_{hh}\) is a unit.
The adjugate formula therefore gives
\[
  C_{hh}^{-1}\in M_{h_h}(A_0).
\]

For every \(i>h\), comparison of the \((i,h)\) blocks then determines
\begin{equation}\label{eq:J-recursion}
  J_{ih}
  =\left(
      B_{ih}
      +\sum_{\ell>h}\xi^{\ell-h}
         B_{i\ell}I_{\ell h}^{\sigma}
      -\sum_{\ell<h}J_{i\ell}C_{\ell h}
    \right)C_{hh}^{-1}.
\end{equation}
This choice makes the \((i,h)\) block of \(JK\) equal to that of
\(BI^\sigma\), and hence eliminates every block below the \(h\)th pivot.
Equations \eqref{eq:C-recursion} and \eqref{eq:J-recursion} therefore
complete the induction.

After the block columns \(0,\ldots,k\) have been eliminated, define the
remaining columns of \(K\) by
\[
  K=J^{-1}BI^\sigma.
\]
The resulting matrices satisfy \(BI^\sigma=JK\), and the first \(k+1\)
block columns of \(K\) have the required normalized upper triangular form.

Finally, the construction is forced at every stage: the \((h,j)\) blocks
with \(h\leq j\leq k\) determine \(C_{hj}\) through
\eqref{eq:C-recursion}, and the blocks below the pivot determine \(J_{ih}\)
through \eqref{eq:J-recursion}.  Hence both \(J\) and \(K\) are unique.
\end{proof}

The factorization defines
\[
  \alpha:\mathfrak I_k\longrightarrow\mathfrak I_k,
  \qquad \alpha(I)=J(I).
\]
For upper factors, use the normalized norm
\[
 \norm{K-K'}_{\mathrm{nor}}
 =\sup_{0\leq i\leq j\leq k}
   \norm{\xi^{-j}(K_{ij}-K'_{ij})}.
\]

\begin{thm}[Block contraction and fixed-point triangularization]
\label{thm:abstract-contraction}
Under \eqref{eq:ordinary-leading-minors}, for all $I,I'\in\mathfrak I_k$,
\begin{align}
  \norm{\alpha(I)-\alpha(I')}
  &\leq\rho\norm{I-I'},\label{eq:alpha-contraction}\\
  \norm{K(I)-K(I')}_{\mathrm{nor}}
  &\leq\rho\norm{I-I'}.
  \label{eq:K-contraction}
\end{align}
Thus $\alpha$ is a strict contraction of the complete metric space
$\mathfrak I_k$ and has a unique fixed point $P$.  If $Q=K(P)$, then
\begin{equation}\label{eq:fixed-triangularization}
  BP^\sigma=PQ,
  \qquad
  P^{-1}BP^\sigma=Q.
\end{equation}
The first $k+1$ block columns of $Q$ are upper block triangular, with
diagonal blocks
\[
  D_i=\xi^iC_{ii},
  \qquad C_{ii}\in\operatorname{GL}_{h_i}(A_0),
  \qquad 0\leq i\leq k.
\]
\end{thm}

\begin{proof}
Write
\[
  BI^\sigma=JK,
  \qquad
  BI'{}^\sigma=J'K'.
\]
For \(0\leq h\leq j\leq k\), let
\(C_{hj}=\xi^{-j}K_{hj}\) and
\(C'_{hj}=\xi^{-j}K'_{hj}\) be the normalized blocks of the two upper
factors.  Put \(\delta=\norm{I-I'}\).  We prove simultaneously, by induction
on \(h\), that
\begin{equation}\label{eq:inductive-contraction-bounds}
  \norm{C_{hj}-C'_{hj}}\leq\rho\delta
  \quad(h\leq j\leq k),
  \qquad
  \norm{J_{ih}-J'_{ih}}\leq\rho\delta
  \quad(i>h).
\end{equation}

Consider first $h=0$.  Subtracting the displayed definitions of $C_{0j}$
and $C'_{0j}$ cancels the terms $B_{0j}$.  Every remaining term contains
\(I_{\ell j}^\sigma-I_{\ell j}'{}^\sigma\), whose norm is at most
\(\delta\), together with the factor \(\xi^{\ell-j}\).  Since
\(\ell>j\), this factor has norm at most \(\rho\).  The ultrametric
inequality therefore gives the first estimate in
\eqref{eq:inductive-contraction-bounds} for \(h=0\).

The same comparison in the displayed formulas defining $J_{i0}$ and
$J'_{i0}$ shows that their parenthesized expressions differ by
an element of norm at most $\rho\delta$.  By \Cref{lem:block-lu}, the pivots
\(C_{00}\) and \(C'_{00}\), as well as their inverses, have entries in
\(A_0\), and hence norm at most one.  Moreover,
\[
  C_{00}^{-1}-C_{00}'{}^{-1}
  =C_{00}^{-1}(C'_{00}-C_{00})C_{00}'{}^{-1},
\]
so the already proved estimate for \(C_{00}-C'_{00}\) gives the same
bound for the difference of their inverses.  Substituting these bounds
into the formulas for \(J_{i0}\) and \(J'_{i0}\) proves the second
estimate for \(h=0\).

Now suppose that the two estimates have been proved for all block columns
of index less than \(h\).  Subtract the two instances of
\eqref{eq:C-recursion}.  Again, the terms \(B_{hj}\) cancel.  The
difference of the first sums has norm at most \(\rho\delta\), for exactly
the same reason as in the initial step.  In the last sum, each difference
\[
  J_{h\ell}C_{\ell j}-J'_{h\ell}C'_{\ell j},
  \qquad \ell<h,
\]
is written as the sum of
\((J_{h\ell}-J'_{h\ell})C_{\ell j}\) and
\(J'_{h\ell}(C_{\ell j}-C'_{\ell j})\).  The factors containing a
difference are bounded by \(\rho\delta\) by the induction hypothesis,
while all remaining factors have norm at most one.  Since the sum over
\(\ell<h\) is finite, the ultrametric inequality proves
\(\norm{C_{hj}-C'_{hj}}\leq\rho\delta\) for every \(j\geq h\).

Let \(A_{ih}\) and \(A'_{ih}\) denote the parenthesized expressions in
\eqref{eq:J-recursion} for \(I\) and \(I'\), respectively.  Comparing
their terms in the same way shows that
\(\norm{A_{ih}-A'_{ih}}\leq\rho\delta\).  The estimate just proved for
the pivots, together with
\[
  C_{hh}^{-1}-C_{hh}'{}^{-1}
  =C_{hh}^{-1}(C'_{hh}-C_{hh})C_{hh}'{}^{-1},
\]
shows that their inverses also differ by at most \(\rho\delta\).
Because \(A'_{ih}\), \(C_{hh}^{-1}\), and \(C_{hh}'{}^{-1}\) all have
norm at most one, comparison of
\(J_{ih}=A_{ih}C_{hh}^{-1}\) and
\(J'_{ih}=A'_{ih}C_{hh}'{}^{-1}\) gives
\(\norm{J_{ih}-J'_{ih}}\leq\rho\delta\).  This completes the induction.

Taking the supremum over the relevant lower blocks gives
\[
  \norm{\alpha(I)-\alpha(I')}
  =\norm{J-J'}
  \leq\rho\delta.
\]
Likewise, since \(K_{hj}=\xi^jC_{hj}\) and
\(K'_{hj}=\xi^jC'_{hj}\), taking the supremum over
\(0\leq h\leq j\leq k\) gives
\[
  \norm{K(I)-K(I')}_{\mathrm{nor}}
  \leq\rho\delta.
\]
Since \(0<\rho<1\), the map \(\alpha\) is a strict contraction of
\(\mathfrak I_k\).  The contraction mapping theorem gives its unique fixed
point $P$.  Substituting $I=P=\alpha(P)$ into the normalized factorization
of \Cref{lem:block-lu} gives \eqref{eq:fixed-triangularization} and the
stated form of the diagonal blocks.
\end{proof}

\subsection{Hodge--Newton filtration and fibrewise eigenvalues}

The following statement is the ordinary Hodge--Newton decomposition.  In
Wan's standard setting it is \cite[Theorem~4.9]{W00a}, where it is attributed
to Dwork; Wan also recalls Katz's more general Hodge--Newton decomposition.
Thus the existence and structure of this filtration are not new.  We include
the proof because the fixed point of the contraction gives a direct
alternative construction of the same filtration and supplies the matrices
needed for the approximation theorem below.

\begin{thm}[Ordinary Hodge--Newton decomposition]
\label{thm:slope-filtration}
Under \eqref{eq:ordinary-leading-minors}, there is an increasing
$\phi$-stable filtration by finite free $A_0$-submodules
\[
  0\subset M_0\subset M_1\subset\cdots\subset M_k\subset M
\]
which is transversal to the basis filtration:
\begin{equation}\label{eq:transversality}
  M=M_i\oplus M(i+1),\qquad0\leq i\leq k.
\end{equation}
Furthermore,
\[
 (M_i/M_{i-1},\phi)\cong(U_i,\xi^i\phi_i),
\]
where $(U_i,\phi_i)$ is a unit root $\sigma$-module of rank $h_i$ whose
Frobenius matrix is $C_{ii}$.  We call $(U_i,\phi_i)$ the unit root
$\sigma$-module coming from the slope-$i$ part of $(M,\phi)$; here
$M_{-1}=0$.
\end{thm}

\begin{proof}
By the discussion following \eqref{eq:I-space-norm}, $P$ is a valid formal
basis change.  Let $\bm f=\bm eP$, and define $M_i$ to be the span of the
first $H_i$ elements of $\bm f$.  Since $Q$ is upper triangular in its first
$k+1$ block columns, the image under Frobenius of every basis vector in
$M_i$ is again in $M_i$.  Thus every $M_i$ is $\phi$-stable.

The leading $H_i\times H_i$ block of $P$ is lower unitriangular, and the
upper-right block from the first $H_i$ rows to the remaining columns is zero.
Consequently, projection along $M(i+1)$ maps the first $H_i$ vectors of
$\bm f$ isomorphically onto the first $H_i$ vectors of $\bm e$.  This proves
\eqref{eq:transversality}.  Finally, the matrix induced by $Q$ on
$M_i/M_{i-1}$ is $D_i=\xi^iC_{ii}$.  Since $C_{ii}$ is invertible over
$A_0$, removing the scalar $\xi^i$ gives the unit root module
$(U_i,\phi_i)$.  This proves all assertions.
\end{proof}

\begin{prop}[Fibrewise slope factors from the diagonal blocks]
\label{prop:fibre-formula}
Let $x$ be $\sigma$-periodic of degree $d$.  For $0\leq i\leq k$, put
\[
  \Phi_{i,x}
  =D_i(x)D_i^\sigma(x)\cdots D_i^{\sigma^{d-1}}(x),
  \qquad
  F_{i,x}(T)=\det(1-T\Phi_{i,x}).
\]
Then the complete slope-$i$ factor of the characteristic series of
$B_{[x]}$ is
\begin{equation}\label{eq:fibre-block-factor}
  F_{i,x}(T)=\det\left(
    1-TD_i(x)D_i^\sigma(x)\cdots D_i^{\sigma^{d-1}}(x)
  \right).
\end{equation}
If $h_i=1$, the corresponding Frobenius eigenvalue is
\begin{equation}\label{eq:fibre-rank-one-root}
  \gamma_{i,x}=\prod_{r=0}^{d-1}D_i^{\sigma^r}(x).
\end{equation}
\end{prop}

\begin{proof}
Apply $\sigma^r$ to \eqref{eq:fixed-triangularization} and multiply the
resulting identities for $0\leq r<d$.  All intermediate $P^{\sigma^r}$
factors cancel, giving
\[
  B B^\sigma\cdots B^{\sigma^{d-1}}P^{\sigma^d}
  =P Q Q^\sigma\cdots Q^{\sigma^{d-1}}.
\]
After specialization at $x$, condition \eqref{eq:periodic-specialization}
identifies $P^{\sigma^d}(x)$ with $P(x)$.  Hence $B_{[x]}$ is similar to
the product of the specialized upper triangular matrices on the right.  The
$i$th diagonal block of a product of block upper triangular matrices is the
ordered product of their $i$th diagonal blocks, which gives
\eqref{eq:fibre-block-factor}.  When $h_i=1$, its determinant is the scalar
product \eqref{eq:fibre-rank-one-root}.
\end{proof}

\subsection{Schur approximation of the graded Frobenius blocks}

The fixed point realizes the slope filtration but need not be explicitly
known.  We now express the finite approximations to its graded Frobenius
blocks directly in terms of the iterated matrices of $B$.

For a block matrix $A$, write
\[
  A_{<i,<i}=(A_{rs})_{0\leq r,s<i},
\]
and use analogous notation for the corresponding block row and block
column.  Whenever $A_{<i,<i}$ is invertible over $A_0[\xi^{-1}]$, define
\[
  \SC_i(A)
  =A_{ii}-A_{i,<i}A_{<i,<i}^{-1}A_{<i,i},
\]
with the convention $\SC_0(A)=A_{00}$.  Thus $\SC_i(A)$ is the last pivot
obtained by block Gaussian elimination on the finite leading truncation
through block $i$.  Write
\[
  \widehat\Delta_i(A)
  :=\det\SC_i(A)
  =\frac{\Delta_i(A)}{\Delta_{i-1}(A)},
  \qquad \Delta_{-1}(A)=1.
\]
More explicitly, if $r$ and $s$ are scalar row and column indices in the
$i$th block, then
\[
  \bigl(\SC_i(A)\bigr)_{rs}
  =\frac{
    \det\begin{pmatrix}
      A_{<i,<i}&A_{<i,s}\\
      A_{r,<i}&A_{rs}
    \end{pmatrix}}
   {\det A_{<i,<i}}.
\]
Thus the Schur complement is determined by bordered minors of the finite
leading truncation of $A$.

Start with $I_0=1$ and define
\begin{equation}\label{eq:contraction-iterates}
  I_{a+1}=\alpha(I_a),
  \qquad
  BI_a^\sigma=I_{a+1}K_{a+1}.
\end{equation}
Let $D_i^{(a)}=(K_a)_{ii}$, and put
\[
  B_0=1,
  \qquad
  B_a=BB^\sigma\cdots B^{\sigma^{a-1}}
  \quad(a\geq1).
\]

\begin{thm}[Schur formula for the finite diagonal blocks]
\label{thm:abstract-schur-formula}
For every $a\geq1$ and $0\leq i\leq k$, the Schur complements appearing
below are defined, and
\[
  \xi^{-ai}\SC_i(B_a)\in\operatorname{GL}_{h_i}(A_0).
\]
\begin{equation}\label{eq:schur-quotient-block}
  D_i^{(a)}
  =\SC_i(B_a)\bigl(\SC_i(B_{a-1})^\sigma\bigr)^{-1}.
\end{equation}
If $h_i=1$, then
\[
  D_i^{(a)}
  =\frac{\widehat\Delta_i(B_a)}
         {\widehat\Delta_i(B_{a-1})^\sigma}.
\]
\end{thm}

\begin{proof}
Iteration of \eqref{eq:contraction-iterates} gives
\begin{equation}\label{eq:iterated-factorization}
  B_a=I_aK_aK_{a-1}^\sigma\cdots K_1^{\sigma^{a-1}}.
\end{equation}
Writing the right side of \eqref{eq:iterated-factorization} as a product,
put
\[
  R_a=K_{a-1}^\sigma\cdots K_1^{\sigma^{a-1}},
  \qquad R_1=1.
\]
Then
\[
  B_a=I_aK_aR_a,
  \qquad
  B_{a-1}^\sigma=I_{a-1}^\sigma R_a.
\]
For $i\leq k$, the leading truncations of $K_a$ and $R_a$ through block
$i$ are upper block triangular and invertible over $A_0[\xi^{-1}]$; the
corresponding truncations of $I_a$ and $I_{a-1}^\sigma$ are lower block
unitriangular.  Hence all the required Schur complements are defined.

Left multiplication by a lower block unitriangular matrix does not alter
a Schur pivot, and the $i$th Schur pivot of an upper block triangular
matrix is its $i$th diagonal block.  It follows that
\[
  \SC_i(B_a)=D_i^{(a)}(R_a)_{ii},
  \qquad
  \SC_i(B_{a-1}^\sigma)=(R_a)_{ii}.
\]
Since Schur complementation commutes with $\sigma$,
\[
  \SC_i(B_{a-1}^\sigma)=\SC_i(B_{a-1})^\sigma.
\]
Moreover,
\[
  (R_a)_{ii}
  =(D_i^{(a-1)})^\sigma
   (D_i^{(a-2)})^{\sigma^2}\cdots
   (D_i^{(1)})^{\sigma^{a-1}}.
\]
Since $\xi^{-i}D_i^{(r)}\in\operatorname{GL}_{h_i}(A_0)$ for every
$r\geq1$, the preceding two displays give
\[
  \xi^{-ai}\SC_i(B_a)\in\operatorname{GL}_{h_i}(A_0).
\]
They also prove \eqref{eq:schur-quotient-block}; together with
the identity $\SC_i(B_a)=\widehat\Delta_i(B_a)$ in rank one, this proves
all the assertions.
\end{proof}

\begin{thm}[Approximation of the graded Frobenius blocks]
\label{thm:abstract-approximation}
For every $a\geq1$ and $0\leq i\leq k$,
\begin{equation}\label{eq:iterate-congruence}
  \xi^{-i}(D_i^{(a)}-D_i)\in\xi^aM_{h_i}(A_0),
\end{equation}
and hence
\begin{equation}\label{eq:iterate-error}
  \norm{\xi^{-i}(D_i^{(a)}-D_i)}\leq\rho^a.
\end{equation}
Combining this estimate with \eqref{eq:schur-quotient-block} gives
\[
  \norm{\xi^{-i}\left(
    \SC_i(B_a)\bigl(\SC_i(B_{a-1})^\sigma\bigr)^{-1}-D_i
  \right)}\leq\rho^a.
\]
If $h_i=1$, this becomes
\[
  \abs{\xi^{-i}\left(
   \frac{\widehat\Delta_i(B_a)}
        {\widehat\Delta_i(B_{a-1})^\sigma}-D_i\right)}
  \leq\rho^a.
\]
\end{thm}

\begin{proof}
Iterating \eqref{eq:alpha-contraction} and using
$\norm{I_0-P}\leq1$ gives
\[
  \norm{I_{a-1}-P}\leq\rho^{a-1}.
\]
The normalized diagonal block of $K_a$ is obtained by applying $K$ to
$I_{a-1}$, whereas $\xi^{-i}D_i$ is obtained by applying $K$ to $P$.
The estimate \eqref{eq:iterate-error} therefore follows from
\eqref{eq:K-contraction}.

The proof of \Cref{thm:abstract-contraction}, with powers of $\xi$ in place
of norm bounds, also gives the ideal-theoretic refinement
\[
  I\equiv I'\pmod{\xi^m}
  \quad\Longrightarrow\quad
  \alpha(I)\equiv\alpha(I')\pmod{\xi^{m+1}},
\]
and the normalized blocks of $K(I)$ and $K(I')$ are congruent modulo
$\xi^{m+1}$.  Indeed, in \eqref{eq:C-recursion} and
\eqref{eq:J-recursion}, every term containing $I-I'$ acquires a factor
$\xi^{\ell-j}$ with $\ell>j$; differences of inverse pivots satisfy the
same congruence by the usual inverse-difference identity.  Applying this
observation to $I_0$ and $P$ proves \eqref{eq:iterate-congruence} and the
corresponding assertion for $I_a-P$.
\end{proof}

Thus the Schur-complement quotient in \eqref{eq:schur-quotient-block}
converges to $D_i$ with normalized error at most $\rho^a$.

There is a direct consequence for the one-step matrix determinant in the
chosen basis.  Write
\[
 \begin{aligned}
  G_i^{(a)}(T)
   &=\det(1-TD_i^{(a)})
     =1+\sum_{m=1}^{h_i}b_{i,m}^{(a)}T^m,\\
  G_i(T)
   &=\det(1-TD_i)
     =1+\sum_{m=1}^{h_i}b_{i,m}T^m.
 \end{aligned}
\]
Since the coefficient of $T^m$ is homogeneous of degree $m$ in the matrix
entries, \eqref{eq:iterate-congruence} gives
\[
  \xi^{-im}\bigl(b_{i,m}^{(a)}-b_{i,m}\bigr)\in\xi^aA_0,
  \qquad 1\leq m\leq h_i.
\]
Equivalently,
\[
  \det(1-T\xi^{-i}D_i^{(a)})
  \equiv\det(1-T\xi^{-i}D_i)\pmod{\xi^aA_0[T]}.
\]
Thus $G_i^{(a)}(T)$ converges coefficientwise to $G_i(T)$.  This is a
matrix-level statement in the chosen basis; the intrinsic fibrewise slope
polynomial is obtained from the periodic Frobenius product considered next.

Taking determinants in \eqref{eq:schur-quotient-block} gives
\begin{equation}\label{eq:minor-ratio-block}
  \det D_i^{(a)}
  =\frac{\widehat\Delta_i(B_a)}
         {\widehat\Delta_i(B_{a-1})^\sigma}.
\end{equation}
Formula \eqref{eq:minor-ratio-block} determines only the top
degree coefficient in higher rank.  The complete polynomial
$G_i^{(a)}(T)$ requires the full matrix-valued Schur quotient in
\eqref{eq:schur-quotient-block}, whose entries are bordered-minor ratios.

\begin{thm}[Fibrewise characteristic-polynomial approximation]
\label{thm:periodic-schur-approximation}
Let $x$ be $\sigma$-periodic of degree $d$, and use
$\Phi_{i,x}$ and $F_{i,x}(T)$ from \Cref{prop:fibre-formula}.
For $a\geq1$, define the ordered finite approximation at consecutive
contraction levels
\begin{equation}\label{eq:consecutive-fibre-approximant}
  \Phi_{i,x}^{(a)}
  =D_i^{(a+d-1)}(x)
   (D_i^{(a+d-2)})^\sigma(x)\cdots
   (D_i^{(a)})^{\sigma^{d-1}}(x).
\end{equation}
Then successive cancellation of the intermediate Schur factors gives
\begin{equation}\label{eq:periodic-schur-cancellation}
  \Phi_{i,x}^{(a)}
  =\SC_i(B_{a+d-1})(x)\SC_i(B_{a-1})(x)^{-1},
\end{equation}
and
\begin{equation}\label{eq:fibre-error}
  \norm{\xi^{-id}
  \bigl(\Phi_{i,x}^{(a)}-\Phi_{i,x}\bigr)}
  \leq\rho^a,
\end{equation}
uniformly in $x$.

Put
\[
  F_{i,x}^{(a)}(T)=\det(1-T\Phi_{i,x}^{(a)}),
\]
and, for $G(T)=\sum_{m=0}^{h_i}g_mT^m\in K_x[T]$, define
\[
  \norm{G}_{i,d}
  :=\max_{0\leq m\leq h_i}\rho^{-idm}\abs{g_m}
  =\norm{G(\xi^{-id}T)}_{\mathrm G}.
\]
Then
\begin{equation}\label{eq:periodic-polynomial-error}
  \norm{F_{i,x}^{(a)}-F_{i,x}}_{i,d}\leq\rho^a.
\end{equation}
\end{thm}

\begin{proof}
Apply \eqref{eq:schur-quotient-block} to the successive factors of
\eqref{eq:consecutive-fibre-approximant}.  For this computation only, put
$E_b=\SC_i(B_b)$.  Before specialization their product is
\[
 \begin{aligned}
  &E_{a+d-1}(E_{a+d-2}^\sigma)^{-1}
   E_{a+d-2}^\sigma(E_{a+d-3}^{\sigma^2})^{-1}\cdots\\
  &\hspace{32mm}\cdots
   E_a^{\sigma^{d-1}}(E_{a-1}^{\sigma^d})^{-1}
   =E_{a+d-1}(E_{a-1}^{\sigma^d})^{-1}.
 \end{aligned}
\]
Since $x\circ\sigma^d=x$, specialization proves
\eqref{eq:periodic-schur-cancellation}.

For $0\leq r<d$, \eqref{eq:iterate-congruence} shows that the normalized
$r$th factor in \eqref{eq:consecutive-fibre-approximant} differs from
$\xi^{-i}D_i^{\sigma^r}(x)$ by an element of
$\xi^{a+d-1-r}M_{h_i}(\mathcal O_{K_x})$, hence by an element of
$\xi^aM_{h_i}(\mathcal O_{K_x})$.  All normalized factors have norm at
most one.  Expanding the difference of the two ordered products proves
\eqref{eq:fibre-error}; the underlying $\xi$-adic congruence makes the
bound uniform in $x$.

After dividing the two matrices by $\xi^{id}$, their $m$th characteristic
coefficients are integral homogeneous polynomials of degree $m$ in the
matrix entries.  Thus, if their characteristic coefficients are denoted by
$c_{i,m}^{(a)}(x)$ and $c_{i,m}(x)$, the same congruence gives
\[
  \abs{\xi^{-idm}
  \bigl(c_{i,m}^{(a)}(x)-c_{i,m}(x)\bigr)}\leq\rho^a,
  \qquad 1\leq m\leq h_i.
\]
This is exactly \eqref{eq:periodic-polynomial-error}.
\end{proof}

\begin{cor}[Rank-one fibre formula]
\label{cor:rank-one-fibre-formula}
If $h_i=1$, write $\gamma_{i,x}=\Phi_{i,x}$ and put
\[
 \begin{aligned}
  \gamma_{i,x}^{(a)}
  &:=\frac{\widehat\Delta_i(B_{a+d-1})(x)}
           {\widehat\Delta_i(B_{a-1})(x)}\\
  &=\frac{
      \Delta_i(B_{a+d-1})\Delta_{i-1}(B_{a-1})
     }{
      \Delta_{i-1}(B_{a+d-1})\Delta_i(B_{a-1})
     }(x).
 \end{aligned}
\]
Then
\[
  \abs{\xi^{-id}
  \bigl(\gamma_{i,x}^{(a)}-\gamma_{i,x}\bigr)}
  \leq\rho^a.
\]
\end{cor}

\begin{proof}
In rank one, \eqref{eq:periodic-schur-cancellation} is the displayed
leading-minor formula, and \eqref{eq:fibre-error} gives its error estimate.
\end{proof}

In particular, \eqref{eq:periodic-polynomial-error} gives the corresponding
coefficientwise Schur-complement limit formula for $F_{i,x}(T)$.  The top
coefficient of the finite approximant satisfies
\begin{equation}\label{eq:periodic-top-coefficient}
  [T^{h_i}]F_{i,x}^{(a)}(T)
  =(-1)^{h_i}
   \frac{\widehat\Delta_i(B_{a+d-1})(x)}
        {\widehat\Delta_i(B_{a-1})(x)}.
\end{equation}
Formula \eqref{eq:periodic-top-coefficient} shows that, for $h_i>1$,
leading principal minors determine only the top coefficient; the full
Schur complements, whose entries are bordered-minor ratios, give the
complete finite polynomial.  The construction does not canonically label
the individual eigenvalues of the corresponding matrix.

The abstract results stated in the Introduction are now completely proved:
\Cref{lem:block-lu} gives the normalized block elimination,
\Cref{thm:abstract-contraction} gives the contraction and its unique fixed
point, \Cref{thm:slope-filtration} identifies the Hodge--Newton filtration,
\Cref{prop:fibre-formula} identifies the fibrewise slope factors,
\Cref{thm:abstract-schur-formula,thm:abstract-approximation} give the exact
finite Schur quotient and its convergence, and
\Cref{thm:periodic-schur-approximation,cor:rank-one-fibre-formula} give the
characteristic-polynomial and rank-one approximations.

The construction above is purely $\sigma$-linear and applies to both
finite-rank and countable-rank nuclear $\sigma$-modules.  In the
countable-rank case, for each prescribed $k$, the nontrivial basis changes
are supported in the first $k+1$ finite block columns, while the nuclear
tail enters through convergent sums.  The construction therefore gives the
corresponding slope filtration and approximations with geometrically
decreasing error in both the finite-rank and countable-rank settings,
without using a connection.

The applications in \Cref{sec:dwork,sec:application} are carried out at the
chain level, on the orthonormal Dwork module $\mathcal N$ and its dual.
Since the contraction theorem also applies in finite rank, it suggests
possible applications at the cohomology level.  We do not pursue a general
study of such applications here.  Instead,
\Cref{sec:differential-interpretation} considers a simple one-parameter
finite-rank setting equipped with a connection compatible with Frobenius,
illustrating how the connection can provide an additional description of
the contraction fixed point and lead to exact formulas for the graded
Frobenius.

\section{The relative Dwork trace formula and duality}\label{sec:dwork}

We now attach a relative Frobenius operator to the family
\eqref{eq:laurent-polynomial}.  The construction is recalled in enough detail
to identify its matrix, its specialization at every closed point, and the
dual nuclear $\sigma$-module to which \Cref{sec:sigma-modules} will be
applied.

\subsection{The relative Dwork operator}

Let $S(f)\subset\bN^m\times\bZ^n$ be the support of $f$, and let
$\Delta=\Delta(f)$ be the convex hull in $\bR^n$ of the origin and the
projection of $S(f)$ to the second factor.  Denote by $C(\Delta)$ the cone
generated by $\Delta$ and put
\[
  M(\Delta)=C(\Delta)\cap\bZ^n.
\]
For $v\in C(\Delta)$, its weight is
\[
  w(v)=\inf\{c\geq0:v\in c\Delta\}.
\]
There is a least integer $D\geq1$ such that
$w(M(\Delta))\subset D^{-1}\bZ_{\geq0}$. Since \(\Delta\) is a polytope
containing the origin, for every \(c\geq0\) one has
\[
\{v\in M(\Delta):w(v)\leq c\} =c\Delta\cap\mathbb Z^n,
\]
and this set is finite.

Let $\bQ_q/\bQ_p$ be the unramified extension with residue field $\bF_q$,
and let $\tau$ be its Frobenius automorphism.  We use the same letter for
the extension of $\tau$ to coefficient fields below.  Recall the
Artin--Hasse exponential
\[
  E(z)=\exp\left(\sum_{r\geq0}\frac{z^{p^r}}{p^r}\right)
  \in1+z\bZ_p[[z]],
\]
where the exponential on the right is understood as a formal power series.
Choose $\pi$ such that $E(\pi)=\zeta_p$, and choose
$\widetilde\pi$ with $\widetilde\pi^D=\pi$.  We extend $\tau$ by
$\tau(\widetilde\pi)=\widetilde\pi$, and write
\[
  K_0=\bQ_q(\widetilde\pi),
  \qquad \theta(z)=E(\pi z).
\]

Let the symbols $a_{u,v}$ also denote the Teichm\"uller lifts of the
coefficients in \eqref{eq:laurent-polynomial}.  Define
\begin{align}
  F(\Lambda,X)
    &=\prod_{(u,v)\in S(f)}
        \theta(a_{u,v}\Lambda^uX^v),\label{eq:splitting-function}\\
  F_a(\Lambda,X)
    &=\prod_{r=0}^{a-1}
       F^{\tau^r}(\Lambda^{p^r},X^{p^r}).
       \notag
\end{align}

Let $\cO=K_0\langle\Lambda_1,\ldots,\Lambda_m\rangle$ be the algebra of
strictly convergent power series, with its Gauss norm, and let
$\cO^\circ$ be its unit ball.  The endomorphism used throughout is
\[
  \sigma(A)(\Lambda)=A^\tau(\Lambda^p).
\]
It is isometric.  The element
\[
  \xi=\widetilde\pi^{p-1}
\]
is fixed by $\sigma$ and satisfies $\ord(\xi)=1/D$.

For real numbers $b>0$ and $c$, let $\cL(b,c)$ be the space of series
\[
  \sum_{v\in M(\Delta)} c_v(\Lambda)X^v,
  \qquad
  c_v(\Lambda)\in\cO,
  \qquad
  \ord(c_v)\geq b\,w(v)+c,
\]
where $\ord(c_v)$ denotes the Gauss valuation.  Put
$\cL(b)=\bigcup_c\cL(b,c)$.  For
$t\in D^{-1}\bZ$, write $\pi^t=\widetilde\pi^{Dt}$.  The orthonormal
Dwork module is
\[
  \cN=
  \left\{
    \sum_{v\in M(\Delta)}c_v(\Lambda)
       \pi^{w(v)}X^v:
    c_v\in\cO,\ c_v\longrightarrow0
  \right\}.
\]
Its orthonormal basis is $e_v=\pi^{w(v)}X^v$.  Let
\[
  \Psi(\sum c_vX^v)=\sum c_{pv}X^v.
\]

First observe that, for every $a\geq1$,
\begin{equation}\label{eq:F-growth}
  F_a(\Lambda,X)\in
  \cL\left(\frac{1}{p^{a-1}(p-1)},0\right).
\end{equation}
Consequently,
\[
  \mathcal T_a=\Psi^a\circ F_a:\cN\longrightarrow\cN
\]
is a completely continuous $\cO$-linear operator.  Indeed, write
$\theta(z)=\sum_{r\geq0}\theta_rz^r$.  The Artin--Hasse estimate is
\begin{equation}\label{eq:theta-coefficient-bound}
  \ord(\theta_r)\geq\frac{r}{p-1}.
\end{equation}
Consider a monomial obtained by choosing powers $r_{u,v}$ from the factors in
\eqref{eq:splitting-function}.  Its $X$-exponent is
$z=\sum r_{u,v}v$.  Since every $v$ in the support lies in $\Delta$ and the
weight is subadditive,
\[
  w(z)\leq\sum_{u,v}r_{u,v}.
\]
Together with \eqref{eq:theta-coefficient-bound}, this shows that the
coefficient of $X^z$ in $F$ has valuation at least $w(z)/(p-1)$; hence
$F\in\cL(1/(p-1),0)$.

Replacing $X$ by $X^{p^r}$ multiplies weights by $p^r$ without changing the
valuations of coefficients.  Therefore
\[
 F^{\tau^r}(\Lambda^{p^r},X^{p^r})
 \in\cL\left(\frac{1}{p^r(p-1)},0\right).
\]
The spaces $\cL(b,0)$ are closed under multiplication, with the smaller of
the two $b$'s controlling a product.  The smallest value among
$r=0,\ldots,a-1$ is $1/(p^{a-1}(p-1))$, proving
\eqref{eq:F-growth}.

Multiplication by $F_a$ sends $\cN$ into
$\cL(p/(p^a(p-1)))$, while
$\Psi^a(\cL(b,c))\subseteq\cL(p^ab,c)$.  Thus
\[
 \cN\xrightarrow{F_a}
 \cL\left(\frac{p}{p^a(p-1)}\right)
 \xrightarrow{\Psi^a}
 \cL\left(\frac{p}{p-1}\right)\subset\cN.
\]
The final inclusion gains decay with respect to weight: coefficients of
output basis vectors of weight tending to infinity tend uniformly to zero.
Truncating the output to the finitely many $v$ with $w(v)\leq R$ therefore
gives finite-rank operators converging to $\mathcal T_a$ in operator norm.
This proves that $\mathcal T_a$ is completely continuous.

We write
\[
  D_a(\Lambda,T)=\det(1-T\mathcal T_a\mid\cN).
\]
Write
\[
  F(\Lambda,X)=\sum_{z\in\bN^m,\,r\in M(\Delta)}
    b_r^{(z)}\Lambda^z\pi^{w(r)}X^r,
  \qquad b_r(\Lambda)=\sum_zb_r^{(z)}\Lambda^z.
\]
In the ordered weight basis, define the infinite matrix
\begin{equation}\label{eq:dwork-matrix}
  B(\Lambda)=(B_{u,v}(\Lambda))_{u,v\in M(\Delta)},
  \qquad
  B_{u,v}=
  \begin{cases}
  b_{pv-u}(\Lambda)\pi^{w(pv-u)+w(u)-w(v)},
     &pv-u\in M(\Delta),\\
  0,&\text{otherwise}.
  \end{cases}
\end{equation}
Indeed, multiplying the input vector $e_u=\pi^{w(u)}X^u$ by $F$ and applying
$\Psi$ produces the output vector $e_v$ precisely from the term with exponent
$pv-u$.  The normalization by $e_v$ gives the coefficient in
\eqref{eq:dwork-matrix}.  Thus the row-$v$, column-$u$ entry of $\mathcal T_1$ is
$B_{u,v}$, so $B(\Lambda)^t$ is its matrix.  Factoring $\mathcal T_a$ into $a$
one-step operators and using
$\Psi(h(X^p)g(X))=h(X)\Psi(g(X))$ gives the ordered matrix product below.
Put
\begin{equation}\label{eq:dwork-matrix-iterate}
  B_0=1,
  \qquad
  B_a(\Lambda)=B(\Lambda)B^\sigma(\Lambda)\cdots
               B^{\sigma^{a-1}}(\Lambda)\quad(a\geq1).
\end{equation}
Then $B_a^t$ is the matrix of $\mathcal T_a$ and hence
\begin{equation}\label{eq:determinant-matrix-identity}
  D_a(\Lambda,T)=\det(1-TB_a(\Lambda)).
\end{equation}

\subsection{The relative trace formula}

For a power series $H(T)$ with constant term one, define
\[
  H(T)^{\delta_a}=\frac{H(T)}{H(p^aT)}.
\]
The following form of Dwork's trace formula is convenient because it is
valid directly on the Dwork Banach module and therefore requires no nondegeneracy
hypothesis.

\begin{thm}[Relative Dwork trace formula]\label{thm:relative-trace-formula}
For $\lambda\in\overline{\bF}_q^m$, let
$d=d(\lambda)$.  Evaluating at the Teichm\"uller lift of $\lambda$, one has
\begin{equation}\label{eq:relative-trace-formula}
  D_d(\lambda,T)^{\delta_d^n}
  =L_\lambda(T)^{(-1)^{n-1}}.
\end{equation}
\end{thm}

\begin{proof}
Let $\lambda$ denote its Teichm\"uller lift.  Since
$d=[\bF_q(\lambda):\bF_p]$, one has $\lambda^{p^d}=\lambda$, and the
specialization $\mathcal T_d(\lambda)$ is a completely continuous linear
operator on the specialized Dwork space.

For $r\geq1$, iterating $\mathcal T_d(\lambda)$ $r$ times gives the
$dr$-step operator.  If $x\in(\bF_{p^{dr}}^*)^n$ and the same letter denotes
its Teichm\"uller lift, the splitting-function identity gives
\begin{align*}
  F_{dr}(\lambda,x)
  &=\prod_{j=0}^{dr-1}\prod_{(u,v)\in S(f)}
      \theta\!\left((a_{u,v}\lambda^u x^v)^{p^j}\right) \\
  &=\zeta_p^{
      \operatorname{Tr}_{\bF_{p^{dr}}/\bF_p}
        (f(\lambda,x))}.
\end{align*}
Thus the usual diagonal trace computation yields
\begin{equation}\label{eq:dwork-trace-identity}
  S_{\lambda,r}=(p^{dr}-1)^n
    \operatorname{Tr}\bigl(\mathcal T_d(\lambda)^r\bigr),
  \qquad r\geq1,
\end{equation}
see \cite{Dw60,AS89}.

Since
\[
  \log D_d(\lambda,T)
  =-\sum_{r\geq1}
    \operatorname{Tr}\bigl(\mathcal T_d(\lambda)^r\bigr)
    \frac{T^r}{r},
\]
applying $\delta_d^n$ multiplies the coefficient of $T^r$ by
$(1-p^{dr})^n$.  Equation \eqref{eq:dwork-trace-identity} therefore gives
\[
  \log D_d(\lambda,T)^{\delta_d^n}
  =(-1)^{n-1}\sum_{r\geq1}\frac{S_{\lambda,r}}{r}T^r
  =\log L_\lambda(T)^{(-1)^{n-1}}.
\]
Exponentiating proves \eqref{eq:relative-trace-formula}.
\end{proof}

Fix a closed point $\lambda$ and put $d=d(\lambda)$.  By
\Cref{thm:relative-trace-formula}, the trace formula can be separated slope
by slope.  If $H(T)$ is a
$p$-adic entire or rational function with constant term one, write
$[H(T)]_s$ for the product of its reciprocal-root and reciprocal-pole factors
of normalized slope $s$, where valuations are divided by $d$.  Empty slope
factors are understood to be one.

\begin{prop}[Slope extraction from the trace formula]
\label{prop:slope-extraction}
Let $\epsilon=(-1)^{n-1}$ and, for $s\in\bQ$, put
\[
  D_{\lambda,s}(T)=[D_d(\lambda,T)]_s,
  \qquad
  \cL_{\lambda,s}(T)=[L_\lambda(T)^\epsilon]_s.
\]
Set $D_{\lambda,s}(T)=1$ for $s<0$.  Then, for every $s\geq0$,
\begin{equation}\label{eq:all-slope-trace-factor}
  \cL_{\lambda,s}(T)
  =\prod_{\ell=0}^{n}
    D_{\lambda,s-\ell}(p^{d\ell}T)^{
       (-1)^\ell\binom{n}{\ell}}.
\end{equation}
In particular, when $0\leq s<1$ only the term $\ell=0$ occurs, so the
slope-$s$ factors of $D_d(\lambda,T)$ and $L_\lambda(T)^\epsilon$ agree.
\end{prop}

\begin{proof}
For any $H(T)$ with constant term one, induction on $n$ gives the binomial
identity
\begin{equation}\label{eq:delta-binomial-identity}
  H(T)^{\delta_d^n}
  =\prod_{\ell=0}^{n}
    H(p^{d\ell}T)^{(-1)^\ell\binom{n}{\ell}}.
\end{equation}
Indeed, applying $\delta_d$ once more divides the product by the same product
with $T$ replaced by $p^dT$; Pascal's identity gives the exponents for
$n+1$.

If $\omega$ is a reciprocal root or pole of $H(T)$ of normalized slope $t$,
then $p^{d\ell}\omega$ is the corresponding reciprocal root or pole of
$H(p^{d\ell}T)$ and has normalized slope $t+\ell$.  Therefore the
slope-$s$ part of the $\ell$th factor on the right of
\eqref{eq:delta-binomial-identity} comes exactly from the slope-$(s-\ell)$
part of $H$.  Taking $H=D_d(\lambda,T)$ and using
\eqref{eq:relative-trace-formula} proves
\eqref{eq:all-slope-trace-factor}.
\end{proof}

\subsection{The dual operator}

The transpose matrix in \eqref{eq:dwork-matrix} is natural for complete
continuity, whereas the untransposed matrix is natural for nuclear
$\sigma$-modules.  We make this passage precise.  Define the formal dual
integral module
\[
  \cN_0^*=\left\{
   \sum_{v\in M(\Delta)}c_v^*(\Lambda)
       \pi^{-w(v)}X^{-v}: c_v^*\in\cO^\circ
  \right\}.
\]
This is a Banach $\cO^\circ$-module with the formal basis
$e_v^*=\pi^{-w(v)}X^{-v}$ in the precise sense of
\Cref{eq:formal-basis-module}.  We write
$\cN^*=\cN_0^*\widehat\otimes_{\cO^\circ}\cO$ for the completed scalar
extension.  For
$\eta=\sum c_v^*\pi^{-w(v)}X^{-v}\in\cN_0^*$ and
$g=\sum c_v\pi^{w(v)}X^v\in\cN$, set
\begin{equation}\label{eq:dwork-pairing}
  \langle\eta,g\rangle
  =\sum_vc_v^*c_v=\CT(\eta g).
\end{equation}
The sum converges because $c_v\to0$.

Let $\Phi(X^{-v})=X^{-pv}$, and let $\operatorname{pr}$ denote projection
of a formal Laurent series onto the monomials $X^{-v}$ with
$v\in M(\Delta)$.  Define
\[
  \mathcal T_a^*=\operatorname{pr}\circ F_a\circ\Phi^a:
  \cN_0^*\longrightarrow\cN_0^*.
\]

\begin{prop}[Dwork duality]\label{prop:dwork-duality}
The operators $\mathcal T_a$ and $\mathcal T_a^*$ are adjoint under
\eqref{eq:dwork-pairing}.  In the formal basis
$e_v^*=\pi^{-w(v)}X^{-v}$, the matrix of $\mathcal T_a^*$ is $B_a$.  In
particular,
\begin{equation}\label{eq:dual-determinant}
  \det(1-T\mathcal T_a\mid\cN)
  =\det(1-T\mathcal T_a^*\mid\cN^*).
\end{equation}
\end{prop}

\begin{proof}
For any Laurent series $h$ for which the following constant terms converge,
the monomials discarded by $\operatorname{pr}$ have exponents outside
$-M(\Delta)$ and therefore cannot contribute to the constant term of the
product with $g\in\cN$.  Hence
\[
  \CT\bigl(\operatorname{pr}(h)g\bigr)=\CT(hg).
\]
If $h=\sum h_vX^{-v}$ and $k=\sum k_uX^u$, then
\[
 \CT(\Phi^a(h)k)
 =\sum_vh_vk_{p^av}
 =\CT(h\Psi^a(k)).
\]
Taking $h=\eta$ and $k=F_ag$ gives
\begin{align*}
 \langle \mathcal T_a^*\eta,g\rangle
 &=\CT\bigl(\operatorname{pr}(F_a\Phi^a(\eta))g\bigr)\\
 &=\CT\bigl(\Phi^a(\eta)F_ag\bigr)\\
 &=\CT\bigl(\eta\Psi^a(F_ag)\bigr)
 =\langle\eta,\mathcal T_ag\rangle.
\end{align*}

Write $\mathcal T_ae_u=\sum_vt_{v,u}e_v$.  Adjointness gives
\[
 \langle \mathcal T_a^*e_v^*,e_u\rangle
 =\langle e_v^*,\mathcal T_ae_u\rangle=t_{v,u};
\]
therefore the matrix of $\mathcal T_a^*$ is the transpose of the matrix of
$\mathcal T_a$.  Since the latter is $B_a^t$, the former is $B_a$.  Finally, the coefficient
of $T^r$ in either Fredholm determinant is the convergent sum of the same
order-$r$ principal minors.  This proves \eqref{eq:dual-determinant}.
\end{proof}

The one-step matrix $B$ therefore defines a $\sigma$-semilinear endomorphism
\begin{equation}\label{eq:dual-sigma-module}
  \phi^*(c e_v^*)=\sigma(c)\sum_u B_{u,v}e_u^*
\end{equation}
of the formal module $\cN_0^*$.  In the terminology of
\Cref{def:nuclear-sigma-module}, the vectors $e_v^*$ form its formal basis.
The valuations of the columns of $B$ tend to infinity, so
$\norm{\phi^*(e_v^*)}\to0$ and $(\cN_0^*,\phi^*)$ is a nuclear
$\sigma$-module.  At a point
$\lambda$ of period $d$, the matrix of its local Frobenius is precisely
$B_d(\lambda)$.  Thus \Cref{prop:dwork-duality} identifies the Fredholm
determinant of the local Frobenius on this abstract module with the
determinant occurring in the trace formula.

\section{Application to toric exponential sums}\label{sec:application}

We apply the abstract construction to the dual Dwork $\sigma$-module
\eqref{eq:dual-sigma-module}.  Throughout this section the points of
$M(\Delta)$ are ordered by increasing weight.

\subsection{The basis polygon and chain-level ordinarity}

For $r\geq0$, put
\[
  W(r)=\#\{v\in M(\Delta):w(v)=r/D\},
  \qquad
  s_r=\sum_{j=0}^rW(j),
  \qquad
  N_r=\sum_{j=0}^rjW(j).
\]
The polygon $P(\Delta)$ is the lower convex polygon joining $(0,0)$ and
the points
\[
  \left(s_r,\frac{N_r}{D}\right),\qquad r\geq0.
\]
It is the basis polygon of the dual Dwork module.

For $u,v\in M(\Delta)$ such that $pv-u\in M(\Delta)$, the subadditivity
and homogeneity of the weight function give
\[
  w(pv-u)+w(u)-w(v)\geq(p-1)w(v).
\]
Together with \eqref{eq:dwork-matrix}, this shows that, in weight blocks,
\[
  B=(\xi^j\mathcal B_{ij})_{i,j\geq0},
  \qquad
  \mathcal B_{ij}\in
  M_{W(i)\times W(j)}(\cO^\circ).
\]
Thus $B$ has precisely the form \eqref{eq:abstract-block-matrix}, and
$P(\Delta)$ is a lower bound for every specialized Fredholm Newton polygon.

For a closed point $\lambda$, normalize valuations by dividing by
$d(\lambda)$.

\begin{defn}\label{def:ordinary-family}
Let $k\geq0$.  The family $f(\Lambda,X)$ is \emph{ordinary on the chain
level up to slope $k/D$} if, for every
$\lambda\in\overline{\bF}_q^m$, the Newton polygon of
$D_{d(\lambda)}(\lambda,T)$ agrees with $P(\Delta)$ through the vertex
$(s_k,N_k/D)$.
\end{defn}

When $k<D$, this agrees, by \Cref{prop:slope-extraction}, with ordinarity of
$L_\lambda(T)^{(-1)^{n-1}}$ through slope $k/D$.
For $k\geq D$, the trace formula combines chain-level slope factors with
integer shifts of factors of lower slope, so chain-level ordinarity is the
appropriate hypothesis for
applying the contraction theorem.  The next criterion connects this geometric
definition with the hypothesis of the abstract theorem, without any
restriction on $k$.

\begin{prop}[Ordinarity criterion]\label{prop:ordinarity-criterion}
The family is ordinary on the chain level up to slope $k/D$ if and only if,
for every
$0\leq r\leq k$, the leading principal minor
\begin{equation}\label{eq:normalized-dwork-minor}
  \xi^{-N_r}\Delta_r(B)
\end{equation}
belongs to $(\cO^\circ)^\times$.  Equivalently, the dual Dwork
$\sigma$-module satisfies the leading-minor hypothesis up to slope $k$ in
the sense of \Cref{def:block-ordinary}.
\end{prop}

\begin{proof}
Write
\[
  \det(1-TB)=1+\sum_{s\geq1}c_s(\Lambda)T^s.
\]
The Fredholm expansion expresses $(-1)^sc_s$ as the convergent sum of the
principal minors indexed by $s$-element subsets of the weight basis.  A
minor indexed by a set $S$ is divisible by
\[
 \xi^{\sum_{v\in S}Dw(v)},
\]
because the column indexed by $v$ is divisible by
$\xi^{Dw(v)}$.  Among subsets of cardinality $s_r$, the smallest possible
sum of the weight indices $Dw(v)$ is $N_r$, and it is attained only by the
subset consisting of all basis vectors of weights
$0,1/D,\ldots,r/D$.  This subset contributes
the leading principal minor.  Every other subset has total index at least
$N_r+1$.  Consequently,
\[
  (-1)^{s_r}\xi^{-N_r}c_{s_r}(\Lambda)
  \equiv\xi^{-N_r}\Delta_r(B)\pmod{\xi\cO^\circ}.
\]

After specialization at a closed point $\lambda$, the one-step Fredholm
Newton polygon meets its basis polygon at $(s_r,N_r/D)$ precisely when
$\abs{g_r(\lambda)}=1$, where $g_r$ is the normalized minor in
\eqref{eq:normalized-dwork-minor}.  The stability of the first contact points
under Frobenius iteration \cite[Corollary~2.5]{W93} shows that this is
equivalent to the corresponding contact for the normalized Newton polygon
of $B_{d(\lambda)}(\lambda)$.  By
\eqref{eq:determinant-matrix-identity}, these are the Newton polygons in
\Cref{def:ordinary-family}.  Thus chain-level ordinarity through slope
$k/D$ is equivalent to $\abs{g_r(\lambda)}=1$ for every closed point
$\lambda$ and every $0\leq r\leq k$.

It remains to pass from all fibres to the coefficient ring $\cO^\circ$.
The strictly convergent series $g_r\in\cO^\circ$ has only finitely many
coefficients of norm one, so its reduction
$\overline g_r$ is a polynomial over the residue field.  The condition
$\abs{g_r(\lambda)}=1$ at every Teichm\"uller point says that
$\overline g_r$ has no zero over the algebraic closure.  Hence
$\overline g_r$ is a nonzero constant and thus a unit.  The affinoid unit
criterion \cite[\S5.1.3, Proposition~1]{BGR} then gives
$g_r\in(\cO^\circ)^\times$.  Conversely, such a unit specializes to a unit
at every closed point.  This proves both implications.
\end{proof}

\subsection{Slope formulas and finite-level approximations}

For a matrix $A$ indexed by the weight basis, continue to write
$\Delta_r(A)$ for its leading principal minor of order $s_r$, and set
$\Delta_{-1}(A)=1$.  Recall the iterated Dwork matrices $B_a$ from
\eqref{eq:dwork-matrix-iterate}.

Fix $k\geq0$ and assume throughout this subsection that
$f(\Lambda,X)$ is ordinary on the chain level up to slope $k/D$.  Apply
the contraction to the dual Dwork matrix $B$.  If $P$ is its fixed point,
write $BP^\sigma=PQ$ and put
\[
  \mathscr D_u=Q_{uu},
  \qquad 0\leq u\leq k.
\]
Thus $\mathscr D_u\in M_{W(u)}(\xi^u\cO^\circ)$ is the Frobenius block on
the limiting slope quotient.

For a closed point $\lambda$ of period $d=d(\lambda)$ and
$0\leq u\leq k$, define
\begin{align}
  \Phi_{u,\lambda}
  &=\mathscr D_u(\lambda)
    \mathscr D_u^\sigma(\lambda)\cdots
    \mathscr D_u^{\sigma^{d-1}}(\lambda),
    \label{eq:local-slope-block}\\
  F_{u,\lambda}(T)
  &=\det(1-T\Phi_{u,\lambda}).
    \notag
\end{align}
\begin{prop}[Slope factors from the limiting triangularization]
\label{prop:dwork-slope-factor}
For every closed point $\lambda$ and $0\leq u\leq k$, the polynomial
$F_{u,\lambda}(T)$ is the complete slope-$u/D$ factor of
$D_d(\lambda,T)$.
\end{prop}

\begin{proof}
By
\Cref{prop:ordinarity-criterion,thm:slope-filtration,prop:fibre-formula}
and Dwork duality \eqref{eq:dual-determinant},
$F_{u,\lambda}(T)$ is precisely the characteristic polynomial of the
corresponding complete slope block.
\end{proof}

Set $F_{u,\lambda}(T)=1$ for $u<0$.  Combining
\Cref{prop:dwork-slope-factor,prop:slope-extraction} gives the following
direct consequences of the Dwork trace formula.  If
$\epsilon=(-1)^{n-1}$, then, for every $0\leq r\leq k$,
\begin{equation}\label{eq:all-slope-L-factor}
  [L_\lambda(T)^\epsilon]_{r/D}
  =\prod_{\ell=0}^{\min(n,\lfloor r/D\rfloor)}
    F_{r-\ell D,\lambda}(p^{d\ell}T)^{
       (-1)^\ell\binom{n}{\ell}}.
\end{equation}
In particular, if $0\leq r\leq k$ and $r<D$, then
\begin{equation}\label{eq:slope-below-one-L-factor}
  [L_\lambda(T)^\epsilon]_{r/D}=F_{r,\lambda}(T).
\end{equation}

We now approximate these graded blocks without using the unknown fixed
point.  By \Cref{thm:abstract-schur-formula}, for every $a\geq1$ and
$0\leq u\leq k$, the Schur pivot $\SC_u(B_a)$ is defined and
\[
  \xi^{-au}\SC_u(B_a)
  \in\operatorname{GL}_{W(u)}(\cO^\circ).
\]
The same assertion for $a=0$ follows immediately from $B_0=1$.

For $a\geq1$, let $\mathscr D_u^{(a)}$ be the $u$th diagonal block
produced by the $a$th contraction step.

\begin{thm}[Schur formula for the finite Dwork blocks]
\label{thm:application-schur-formula}
For every $a\geq1$ and $0\leq u\leq k$,
\begin{equation}\label{eq:application-schur-block}
  \mathscr D_u^{(a)}
  =\SC_u(B_a)\bigl(\SC_u(B_{a-1})^\sigma\bigr)^{-1}.
\end{equation}
If $W(u)=1$, this specializes to
\[
  \mathscr D_u^{(a)}
  =\frac{\widehat\Delta_u(B_a)}
         {\widehat\Delta_u(B_{a-1})^\sigma}.
\]
\end{thm}

\begin{proof}
This is \Cref{thm:abstract-schur-formula} for the Dwork matrix.  The final
formula follows because a rank-one Schur complement is the scalar
$\widehat\Delta_u$.
\end{proof}

\begin{thm}[Approximation of the Dwork graded blocks]
\label{thm:application-block-approximation}
For every $a\geq1$ and $0\leq u\leq k$,
\begin{equation}\label{eq:application-block-error}
  \norm{\xi^{-u}(\mathscr D_u^{(a)}-\mathscr D_u)}
  \leq\abs{\xi}^a.
\end{equation}
Consequently,
\[
  \norm{\xi^{-u}\left(
    \SC_u(B_a)\bigl(\SC_u(B_{a-1})^\sigma\bigr)^{-1}
    -\mathscr D_u\right)}
  \leq\abs{\xi}^a.
\]
\end{thm}

\begin{proof}
The first estimate is \Cref{thm:abstract-approximation} in the Dwork
setting; the second follows from \eqref{eq:application-schur-block}.
\end{proof}

For a closed point $\lambda$ of period $d=d(\lambda)$, define
\begin{equation}\label{eq:application-consecutive-approximant}
 \begin{aligned}
  \Phi_{u,\lambda}^{(a)}
  ={}&\mathscr D_u^{(a+d-1)}(\lambda)
      (\mathscr D_u^{(a+d-2)})^\sigma(\lambda)\cdots\\
     &\hspace{28mm}\cdots
      (\mathscr D_u^{(a)})^{\sigma^{d-1}}(\lambda),
  \qquad
  F_{u,\lambda}^{(a)}(T)
  =\det(1-T\Phi_{u,\lambda}^{(a)}).
 \end{aligned}
\end{equation}
For $G(T)=\sum_{m=0}^{W(u)}g_mT^m$ over the specialization field, put
\[
  \norm{G}_{u,d}
  :=\max_{0\leq m\leq W(u)}
       \abs{\xi}^{-udm}\abs{g_m}
  =\norm{G(\xi^{-ud}T)}_{\mathrm G}.
\]

\begin{thm}[Finite-iterate approximation of the Dwork slope polynomials]
\label{thm:application-fibre}
For every $a\geq1$, every closed point $\lambda$ of period $d$, and
$0\leq u\leq k$, successive cancellation of the intermediate Schur
factors in \eqref{eq:application-consecutive-approximant} gives
\begin{equation}\label{eq:application-periodic-approximant}
  \Phi_{u,\lambda}^{(a)}
  =\SC_u(B_{a+d-1})(\lambda)
   \SC_u(B_{a-1})(\lambda)^{-1}.
\end{equation}
Moreover,
\begin{equation}\label{eq:application-fibre-block-error}
  \norm{\xi^{-ud}
    (\Phi_{u,\lambda}^{(a)}-\Phi_{u,\lambda})}
  \leq\abs{\xi}^a,
\end{equation}
and
\begin{equation}\label{eq:application-slope-factor-error}
  \norm{F_{u,\lambda}^{(a)}-F_{u,\lambda}}_{u,d}
  \leq\abs{\xi}^a.
\end{equation}
All three assertions are uniform in $\lambda$.
\end{thm}

\begin{proof}
Apply \Cref{thm:periodic-schur-approximation} to the Dwork matrix $B$,
with $i=u$ and $x$ the Teichm\"uller specialization at $\lambda$.
Its ordered finite approximation is precisely
\eqref{eq:application-consecutive-approximant}; the three conclusions are
\eqref{eq:periodic-schur-cancellation}, \eqref{eq:fibre-error}, and
\eqref{eq:periodic-polynomial-error}, respectively.
\end{proof}

\begin{cor}[Rank-one leading-minor formula]
\label{cor:application-rank-one}
Suppose $W(u)=1$.  Write $\gamma_{u,\lambda}=\Phi_{u,\lambda}$ and set
\begin{equation}\label{eq:application-rank-one-root}
 \begin{aligned}
  \gamma_{u,\lambda}^{(a)}
  :=\Phi_{u,\lambda}^{(a)}
  &=\frac{\widehat\Delta_u(B_{a+d-1})(\lambda)}
          {\widehat\Delta_u(B_{a-1})(\lambda)}\\
  &=\frac{
      \Delta_u(B_{a+d-1})\Delta_{u-1}(B_{a-1})
     }{
      \Delta_{u-1}(B_{a+d-1})\Delta_u(B_{a-1})
     }(\lambda).
 \end{aligned}
\end{equation}
Then $F_{u,\lambda}(T)=1-\gamma_{u,\lambda}T$ and
\begin{equation}\label{eq:application-rank-one-error}
  \abs{\xi^{-ud}
    (\gamma_{u,\lambda}^{(a)}-\gamma_{u,\lambda})}
  \leq\abs{\xi}^a.
\end{equation}
\end{cor}

\begin{proof}
This is \Cref{cor:rank-one-fibre-formula} for the Dwork matrix.
\end{proof}

Thus \Cref{thm:application-schur-formula,thm:application-block-approximation}
give an explicit Schur-complement formula for the finite block and its
limit $\mathscr D_u$; more precisely,
\eqref{eq:application-block-error} controls the coefficient-ring block.
At a periodic fibre, \eqref{eq:application-fibre-block-error} controls the
full ordered Frobenius product, and \Cref{thm:application-fibre} gives the
coefficientwise limit
$F_{u,\lambda}^{(a)}(T)\to F_{u,\lambda}(T)$ with a uniform quantified
error.  The scalar specialization is recorded separately in
\Cref{cor:application-rank-one}.

The same approximation applies coefficientwise to the one-step matrix
determinant polynomials in the chosen basis:
\[
  \det(1-T\xi^{-u}\mathscr D_u^{(a)})
  \equiv
  \det(1-T\xi^{-u}\mathscr D_u)
  \pmod{\xi^a\cO^\circ[T]}.
\]
Thus $\det(1-T\mathscr D_u^{(a)})$ converges coefficientwise to
$\det(1-T\mathscr D_u)$.  This is a coefficient-ring statement for one
$\sigma$-linear block; the fibrewise slope factor uses the full periodic
product.

Equivalently, if
\[
  F_{u,\lambda}^{(a)}(T)
  =1+\sum_{m=1}^{W(u)}c_{u,m,\lambda}^{(a)}T^m,
  \qquad
  F_{u,\lambda}(T)
  =1+\sum_{m=1}^{W(u)}c_{u,m,\lambda}T^m,
\]
then \eqref{eq:application-slope-factor-error} says
\[
  \abs{\xi^{-udm}
  \bigl(c_{u,m,\lambda}^{(a)}-c_{u,m,\lambda}\bigr)}
  \leq\abs{\xi}^a,
  \qquad 1\leq m\leq W(u).
\]
The normalized finite matrix $\xi^{-ud}\Phi_{u,\lambda}^{(a)}$ is invertible
over the valuation ring of the specialization field; hence every reciprocal
root of $F_{u,\lambda}^{(a)}(T)$ has normalized slope $u/D$.  For
$0\leq u\leq k$ with $u<D$, \eqref{eq:slope-below-one-L-factor} identifies
$F_{u,\lambda}(T)$ directly with the slope-$u/D$ factor of
$L_\lambda(T)^\epsilon$.  For larger $u$, it is the complete slope factor
of the chain-level Fredholm determinant entering the trace formula.

The top coefficient of the finite polynomial is
\[
  \bigl[T^{W(u)}\bigr]F_{u,\lambda}^{(a)}(T)
  =(-1)^{W(u)}
   \frac{\widehat\Delta_u(B_{a+d-1})(\lambda)}
        {\widehat\Delta_u(B_{a-1})(\lambda)}.
\]
For $W(u)>1$, leading principal minors determine only this coefficient;
the remaining coefficients require the full Schur complements in
\eqref{eq:application-periodic-approximant}.  For $W(u)=1$, the displayed
identity is exactly \eqref{eq:application-rank-one-root}.

\medskip
\noindent\emph{Passage to the $L$-function.}
Set $F_{u,\lambda}^{(a)}(T)=1$ for $u<0$.  For $0\leq r\leq k$, put
\begin{equation}\label{eq:application-rational-approximant}
  \mathcal R_{r,\lambda}^{(a)}(T)
  =\prod_{\ell=0}^{\min(n,\lfloor r/D\rfloor)}
    F_{r-\ell D,\lambda}^{(a)}
      (p^{d\ell}T)^{(-1)^\ell\binom n\ell}.
\end{equation}
\begin{thm}[Approximation of the $L$-function slope factors]
\label{thm:application-L-factor}
This converges coefficientwise, uniformly in $\lambda$, to the complete
rational slope factor in \eqref{eq:all-slope-L-factor}.  More precisely, if
\[
  \mathcal R_{r,\lambda}^{(a)}(T)
  =\sum_{m\geq0}\kappa_{r,m,\lambda}^{(a)}T^m,
  \qquad
  [L_\lambda(T)^\epsilon]_{r/D}
  =\sum_{m\geq0}\kappa_{r,m,\lambda}T^m,
\]
as formal power series, then
\begin{equation}\label{eq:application-rational-coefficient-error}
  \abs{\xi^{-rdm}
  \bigl(\kappa_{r,m,\lambda}^{(a)}-\kappa_{r,m,\lambda}\bigr)}
  \leq\abs{\xi}^a,
  \qquad m\geq1.
\end{equation}
\end{thm}

\begin{proof}
For a fixed $r$, substitute the finite polynomials into
\eqref{eq:application-rational-approximant}.  For the coefficient of
$T^m$ in the factor indexed by $\ell$, multiplication by
$p^{d\ell m}$ changes the $\xi$-exponent in the normalization from
$(r-\ell D)dm$ to $rdm$.  Since
\[
  \abs{(p^{d\ell}\xi^{-Dd\ell})^m}=1,
\]
the normalized coefficient estimate
\eqref{eq:application-slope-factor-error} is preserved.  Products and
formal inverses of series with constant term one preserve the same
$\xi^a$-congruence, proving
\eqref{eq:application-rational-coefficient-error}.
\end{proof}

\begin{cor}[Rank-one slope formulas]\label{cor:rank-one-all-slopes}
Fix $0\leq r\leq k$ and suppose that
$W(r-\ell D)=1$ for every
$0\leq\ell\leq\min(n,\lfloor r/D\rfloor)$.  For each such $\ell$, use the
scalars $\gamma_{r-\ell D,\lambda}$ and
$\gamma_{r-\ell D,\lambda}^{(a)}$ defined above, and put
\begin{align*}
  \beta_{r,\ell,\lambda}
  &=p^{d\ell}\gamma_{r-\ell D,\lambda},
  \qquad
  \beta_{r,\ell,\lambda}^{(a)}
  =p^{d\ell}\gamma_{r-\ell D,\lambda}^{(a)}.
\end{align*}
Then every $\beta_{r,\ell,\lambda}$ has normalized slope $r/D$, and
\begin{align}
  [L_\lambda(T)^\epsilon]_{r/D}
  &=\prod_{\ell=0}^{\min(n,\lfloor r/D\rfloor)}
    (1-\beta_{r,\ell,\lambda}T)^{
       (-1)^\ell\binom{n}{\ell}},
       \label{eq:rank-one-all-slope-factor}\\
  [L_\lambda(T)]_{r/D}
  &=\prod_{\ell=0}^{\min(n,\lfloor r/D\rfloor)}
    (1-\beta_{r,\ell,\lambda}T)^{
       \epsilon(-1)^\ell\binom{n}{\ell}}.
       \label{eq:rank-one-all-slope-factor-L}
\end{align}
Moreover, for every $a\geq1$ and
$0\leq\ell\leq\min(n,\lfloor r/D\rfloor)$, the finite-minor
approximations satisfy
\begin{equation}\label{eq:fibre-root-error}
  \abs{\xi^{-rd}\bigl(
    \beta_{r,\ell,\lambda}^{(a)}
    -\beta_{r,\ell,\lambda}\bigr)}
  \leq\abs{\xi}^a.
\end{equation}
After common factors of opposite signs in
\eqref{eq:rank-one-all-slope-factor-L} are cancelled, the remaining
$\beta_{r,\ell,\lambda}$ are precisely the reciprocal zeros and poles of
$L_\lambda(T)$ of slope $r/D$, with the displayed multiplicities.  In
particular, if no two contributions of opposite signs coincide, there is no
cancellation.
\end{cor}

\begin{proof}
For each admissible $\ell$, put $u=r-\ell D$.  Then $W(u)=1$, so the local
block \eqref{eq:local-slope-block} is the scalar $\gamma_{u,\lambda}$.  Its
valuation is $ud/D$.  Therefore
$p^{d\ell}\gamma_{r-\ell D,\lambda}$ has valuation $rd/D$, and
\eqref{eq:rank-one-all-slope-factor} follows from
\eqref{eq:all-slope-L-factor}.  Raising both sides to the power
$\epsilon$ gives \eqref{eq:rank-one-all-slope-factor-L}.

By \eqref{eq:application-rank-one-error},
\[
  \abs{\xi^{-(r-\ell D)d}
  \bigl(\gamma_{r-\ell D,\lambda}^{(a)}
        -\gamma_{r-\ell D,\lambda}\bigr)}
  \leq\abs{\xi}^a.
\]
Because $\ord(p)=1$ and $\ord(\xi)=1/D$, the element
$p^{d\ell}\xi^{-\ell Dd}$ has norm one.  Multiplying the last estimate by
this element proves \eqref{eq:fibre-root-error}.  The final assertion is the
usual cancellation of equal linear factors in a rational function.  Such
cancellation need not occur exactly at a finite approximation level.
\end{proof}

Since $F_{u,\lambda}$ has degree $W(u)$, taking signed degrees in
\eqref{eq:all-slope-L-factor} also gives, for $0\leq r\leq k$,
\[
 \begin{aligned}
  \deg [L_\lambda(T)^\epsilon]_{r/D}
  &=\sum_{\ell=0}^{\min(n,\lfloor r/D\rfloor)}
    (-1)^\ell\binom{n}{\ell}W(r-\ell D),\\
  \deg [L_\lambda(T)]_{r/D}
  &=\epsilon
    \sum_{\ell=0}^{\min(n,\lfloor r/D\rfloor)}
    (-1)^\ell\binom{n}{\ell}W(r-\ell D).
 \end{aligned}
\]
Here degree means numerator degree minus denominator degree.

Formula \eqref{eq:slope-below-one-L-factor} shows that an ordinary slope
factor with $r/D<1$ is governed by a single chain-level Frobenius block.
When that block
has rank one, \Cref{cor:application-rank-one} gives a ratio of leading
principal minors of two finite Dwork iterates, extending the unit-root case
$r=0$ to positive slopes below one; \Cref{cor:rank-one-all-slopes} records
the trace-formula combination when all relevant blocks have rank one.  In
higher rank,
\Cref{thm:application-block-approximation,thm:application-fibre} give,
respectively, the finite-iterate approximation to the graded block and the
complete fibrewise slope polynomial with a quantified error;
\Cref{thm:application-L-factor} then gives the corresponding rational
$L$-function slope factor.
For $r\geq D$, the trace formula combines the characteristic factors attached
to chain-level blocks with nonnegative indices among
$r,r-D,\ldots,r-nD$, with alternating binomial
multiplicities.  No new contraction argument is required: the same
construction applies at each
slope for which chain-level ordinarity holds.

\subsection{The unit-root formula of Adolphson and Sperber}
\label{subsec:as-unit-root-comparison}

We show that our contraction map can recover the formula of Adolphson and Sperber for the unit roots. We use their classical splitting function and write $B_{\mathrm{AS}}$ for the resulting
Frobinius matrix, to distinguish it from the Artin--Hasse realization above.

Consider the universal family
\[
 f_\Lambda(X)=\sum_{a\in A}\Lambda_aX^a,
 \qquad A\subset\mathbb Z^n
\]
with $A$ finite and spanning $\mathbb R^n$.  Choose $\varpi$ with
$\varpi^{p-1}=-p$, corresponding to the chosen additive character, and put
\[
 \theta_{\mathrm{AS}}(t)=\exp\bigl(\varpi(t-t^p)\bigr),
 \qquad
 \mathcal E(\Lambda,X)=\exp\bigl(\varpi f_\Lambda(X)\bigr).
\]
Here $\mathcal E$ is the ordinary exponential generating series.  The
splitting factor satisfies the formal identity
\begin{equation}\label{eq:as-splitting-quotient}
 F_{\mathrm{AS}}(\Lambda,X)
 :=\prod_{a\in A}\theta_{\mathrm{AS}}(\Lambda_aX^a)
 =\frac{\mathcal E(\Lambda,X)}
        {\mathcal E(\Lambda^p,X^p)}.
\end{equation}
Let $C$ be the cone generated by $A$, let $M=C\cap\mathbb Z^n$, and let
$\operatorname{pr}_{-}$ retain the monomials with exponents in $-M$.
Define
\[
 \mathcal G(\Lambda,X)=\operatorname{pr}_{-}\mathcal E(\Lambda,X).
\]
Since $\mathcal E$ is supported in $C$, this projection retains precisely
its terms with exponents in $C\cap(-C)$.  In particular,
$\mathcal G=\mathcal E$ if $C=\mathbb R^n$.
The relative dual Frobenius is
\[
 \phi_{\mathrm{AS}}^*(h)
 =\operatorname{pr}_{-}\!\left(
 F_{\mathrm{AS}}(\Lambda,X)h(\Lambda^p,X^p)\right).
\]

We claim that $\phi_{\mathrm{AS}}^*(\mathcal G)=\mathcal G$.
Both $F_{\mathrm{AS}}$ and $\mathcal E$ are supported in $C$.
If $u\in C$ and $v\in C\setminus(-C)$, then $u+v\notin-C$:
otherwise $-v=u-(u+v)\in C$, a contradiction.
Consequently, a term discarded from $\mathcal E(\Lambda^p,X^p)$ cannot
contribute to the projection after multiplication by $F_{\mathrm{AS}}$.
Using \eqref{eq:as-splitting-quotient}, we obtain
\[
 \begin{aligned}
 \phi_{\mathrm{AS}}^*(\mathcal G)
 &=\operatorname{pr}_{-}\!\left(
     F_{\mathrm{AS}}(\Lambda,X)\mathcal G(\Lambda^p,X^p)\right)\\
 &=\operatorname{pr}_{-}\!\left(
     F_{\mathrm{AS}}(\Lambda,X)\mathcal E(\Lambda^p,X^p)\right)\\
 &=\operatorname{pr}_{-}\mathcal E(\Lambda,X)
 =\mathcal G(\Lambda,X).
 \end{aligned}
\]
These identities are first understood formally in the parameters; each
fixed parameter degree involves only finitely many monomials.

To express this calculation in matrix form, use the weight $w$ associated
with $\operatorname{conv}(A\cup\{0\})$, choose $D\ge1$ with
$Dw(M)\subset\mathbb Z_{\ge0}$, and choose $\eta$ with
$\ord(\eta)=(p-1)/p^2$.  Enlarge a finite coefficient field $L$ to contain
$\varpi$, $\eta$, and the roots needed for $\eta^{w(\nu)}$.
In the weighted dual basis
\[
 e_\nu^{\mathrm{AS}}=\eta^{-w(\nu)}X^{-\nu},\qquad \nu\in M,
\]
write
\[
 \mathcal G=\sum_{\nu\in M}G_\nu(\Lambda)e_\nu^{\mathrm{AS}},
 \qquad G=(G_\nu)_{\nu\in M}.
\]
Let $\sigma$ fix $L$ and send $\Lambda_a$ to $\Lambda_a^p$.
By definition of the Frobenius matrix, the preceding identity becomes
\begin{equation}\label{eq:as-invariant-generating-vector}
 B_{\mathrm{AS}}G^\sigma=G.
\end{equation}
The substitution $X\mapsto X^p$ is already incorporated in
$B_{\mathrm{AS}}$; the superscript $\sigma$ here acts only on the
coefficient functions.

The constant coordinate is
\[
 G_0(\Lambda)=\operatorname{CT}_X\mathcal E(\Lambda,X)
 =\sum_{\substack{(m_a)\in\mathbb Z_{\ge0}^{A}\\
                   \sum_a m_aa=0}}
   \frac{\varpi^{\sum_a m_a}\prod_a\Lambda_a^{m_a}}
        {\prod_a m_a!}.
\]
Thus $G_0$ is the $A$-hypergeometric series with parameter zero,
with its variables scaled by $\varpi$.  The coordinates of $G$ belong
to $R=\mathcal O_L[[\Lambda_a:a\in A]]$, endowed with the supremum norm
on its coefficients.  Since $G_0\equiv1\pmod{\varpi R}$, it is a unit
in $R$.  Put $z=G/G_0$.  Its constant coordinate is one, and
\eqref{eq:as-invariant-generating-vector} gives
\begin{equation}\label{eq:as-normalized-generating-vector}
 B_{\mathrm{AS}}z^\sigma
 =\frac{G_0}{G_0^\sigma}z.
\end{equation}

For completeness, the connection with our contraction is explicit.
Separate the constant coordinate from the positive-weight coordinates
and write
\[
 B_{\mathrm{AS}}=
 \begin{pmatrix}b_{00}&b_{0+}\\ b_{+0}&B_{++}\end{pmatrix}.
\]
For $I(v)=\left(\begin{smallmatrix}1&0\\v&1\end{smallmatrix}\right)$,
the first-column LU factorization gives
\[
 \alpha_0(v)
 =\bigl(b_{+0}+B_{++}v^\sigma\bigr)
    \bigl(b_{00}+b_{0+}v^\sigma\bigr)^{-1}.
\]
The coefficient estimates in \cite[Lemmas~2.2 and~2.8]{AS12} give a
column-nuclear matrix with unit $b_{00}$.  More precisely, with
$\delta=\eta^{(p-1)/D}$, its column indexed by $\nu$ is divisible
by $\delta^{Dw(\nu)}$.  Since $w(\nu)>0$ for $\nu\ne0$,
the constant block has rank one and the slope-zero contraction applies;
no ordinarity assumption on positive slopes is needed.
Equation~\eqref{eq:as-normalized-generating-vector} identifies $z$ with
its unique normalized fixed point and yields the scalar graded block
\[
 D_{0,\mathrm{AS}}(\Lambda)
 =\mathcal F(\Lambda)
 :=\frac{G_0(\Lambda)}{G_0(\Lambda^p)}.
\]

The same contraction preserves the corresponding module over
$R'=\mathcal O_L\langle\Lambda_a:a\in A\rangle$.
Its iterates starting from the constant basis vector therefore converge
over $R'$.  Uniqueness over $R$ identifies this limit with $z$, so that
$z$ and $\mathcal F$ extend to the closed unit polydisc.  This is the brilliant 
analytic-continuation mechanism in
\cite[Propositions~2.11 and~2.15]{AS12}.
For a Teichm\"uller parameter $\widehat\lambda$ of $p$-Frobenius period
$d$, the rank-one fibre formula and the Dwork trace formula give
\[
 u(\lambda)
 =\prod_{j=0}^{d-1}D_{0,\mathrm{AS}}^{\sigma^j}(\widehat\lambda)
 =\prod_{j=0}^{d-1}\mathcal F(\widehat\lambda^{p^j}),
\]
the unit-root formula of \cite[Theorem~1.3]{AS12}.
Here $1-u(\lambda)T$ is the complete slope-zero factor of
$L_\lambda(T)^{(-1)^{n-1}}$ for this universal family.
The values in this product are those of the extended quotient
$\mathcal F$; the individual series $G_0$ need not converge at
$\widehat\lambda$.  In particular, the formal identity
\eqref{eq:as-invariant-generating-vector} cannot in general be specialized
directly at that point.

\section{A one-parameter differential-module illustration}
\label{sec:differential-interpretation}

The preceding sections apply the contraction construction directly to the
orthonormal Dwork module $\mathcal N$ and its dual at the chain level.  There,
Schur complements of consecutive iterated Frobenius matrices give convergent
approximations to the graded Frobenius blocks.  The abstract theorem also
applies to finite-rank $\sigma$-modules, where additional structures may
identify the limiting blocks directly.

This section illustrates how a compatible connection can provide exact
expressions for the graded Frobenius matrices constructed by the contraction
theorem.  We work with a finite-rank $\sigma$-module over
$A_0=\mathcal O_K\langle t\rangle$.  The contraction fixed point $P$, the
transformed Frobenius matrix $Q$, and its diagonal blocks $D_i$ are all
defined over $A_0$.  A normalized solution matrix of the associated
differential equation gives local Schur expressions for these global
matrices.  The resulting expressions, taken as a whole, belong to $A_0$
and therefore determine characteristic factors at periodic specializations
of $A_0$.  In rank one, they yield root formulas in terms of leading
principal minors.

This one-parameter example illustrates a possible application of the
contraction construction when a compatible differential structure is
available, as may occur at the cohomological level.  Our principal
application remains the chain-level construction of the preceding
sections; none of the results below is used in Section~4.

Let $K$ be a complete discretely valued field of characteristic zero, let
$\mathcal O_K$ be its valuation ring, and let
\[
  A_0=\mathcal O_K\langle t\rangle,
  \qquad A=K\langle t\rangle.
\]
Equip $A_0$ with its Gauss norm.  Fix an integer $q>1$ and a
$K$-algebra endomorphism $\sigma:A\to A$ preserving $A_0$ such that
\begin{equation}\label{eq:differential-frobenius-lift}
  \sigma|_K=\operatorname{id}_K,
  \qquad \sigma(t)=t^q.
\end{equation}
For $\partial=d/dt$, one then has
\begin{equation}\label{eq:derivation-frobenius-relation}
  \partial\circ\sigma
  =qt^{q-1}\sigma\circ\partial.
\end{equation}

Let $M^\circ$ be a finite free $A_0$-module and put
$M=M^\circ\otimes_{\mathcal O_K}K$.  A connection in the $t$-direction is
an additive map $\nabla_\partial:M\longrightarrow M$ satisfying the Leibniz
rule
\[
  \nabla_\partial(am)=\partial(a)m+a\nabla_\partial(m)
  \qquad(a\in A,\ m\in M).
\]
Let $\phi:M^\circ\to M^\circ$ be $\sigma$-linear, and assume that its
linearization becomes an isomorphism over $A$.  We say that $\phi$ is
\emph{compatible} with $\nabla_\partial$ if
\begin{equation}\label{eq:connection-frobenius-compatibility}
  \nabla_\partial\circ\phi
  =qt^{q-1}\phi\circ\nabla_\partial.
\end{equation}
This is the one-derivation form of the usual requirement that a Frobenius
structure be a morphism of modules with connection.  In particular, this
compatibility implies that $\phi$ sends horizontal elements to horizontal
elements.

Choose a row basis $\bm e$ of $M^\circ$ and write
\[
  \phi(\bm e)=\bm eB,
  \qquad
  \nabla_\partial(\bm e)=\bm eN,
  \qquad
  B\in M_n(A_0)\cap\operatorname{GL}_n(A),
  \quad N\in M_n(A).
\]
Then \eqref{eq:connection-frobenius-compatibility} is equivalent to
\begin{equation}\label{eq:matrix-connection-compatibility}
  \partial(B)+NB=qt^{q-1}BN^\sigma.
\end{equation}
Assume that the basis is grouped into blocks of sizes
$h_0,\ldots,h_k$ with $\sum_{j=0}^k h_j=n$, and that
\begin{equation}\label{eq:differential-weighted-columns}
  B=(\xi^jB_{ij})_{0\leq i,j\leq k},
  \qquad B_{ij}\in M_{h_i\times h_j}(A_0),
\end{equation}
where $\xi\in\mathcal O_K$ satisfies $0<\abs{\xi}<1$.  Put
\[
  H_i=\sum_{j=0}^i h_j,
  \qquad N_i=\sum_{j=0}^i jh_j,
\]
and let $\Delta_i(V)$ denote the leading principal minor of order $H_i$
of a block matrix $V$, with $\Delta_{-1}(V)=1$.  We impose the
leading-minor hypothesis of Section~2:
\[
  \Delta_i(B)=\xi^{N_i}v_i,
  \qquad v_i\in A_0^\times,
  \qquad 0\leq i\leq k.
\]
Together with \eqref{eq:differential-weighted-columns}, this hypothesis
allows us to apply the finite-rank contraction theorem over $A_0$.
It gives a unique lower block unitriangular fixed point
$P\in\operatorname{GL}_n(A_0)$ and an upper block triangular matrix
$Q\in M_n(A_0)\cap\operatorname{GL}_n(A)$ satisfying
\begin{equation}\label{eq:differential-fixed-triangularization}
  BP^\sigma=PQ.
\end{equation}

Writing $D_i=Q_{ii}$, we have
\[
  \xi^{-i}D_i\in\operatorname{GL}_{h_i}(A_0).
\]
Thus $P$, $Q$, and the graded Frobenius matrices $D_i$ are defined over
$A_0$ before any restriction of the parameter disc.  They can be specialized
at every bounded $\mathcal O_K$-algebra homomorphism from $A_0$ to a
valuation ring.  We use the same symbols for their restrictions below.

We next make precise the local base change used to construct a normalized
solution matrix.  For $r\in\abs{K^\times}$ with $0<r\leq1$, put
\[
  A_r=\left\{\sum_{m\geq0}a_mt^m\in K[[t]]:
                \abs{a_m}r^m\longrightarrow0\right\},
\]
with norm and integral unit ball
\[
  \norm{\textstyle\sum_m a_mt^m}_r=\sup_m\abs{a_m}r^m,
  \qquad A_r^\circ=\{f\in A_r:\norm{f}_r\leq1\}.
\]
Restriction from the closed unit disc gives $A_0\to A_r^\circ$.  The ring
$A_r^\circ$ is complete, its norm topology agrees with its $\xi$-adic
topology, and $\sigma$ preserves it.  If $r<1$, the induced endomorphism need
not be isometric, but it is nonexpanding:
\begin{equation}\label{eq:sigma-nonexpanding-small-disc}
  \norm{f^\sigma}_r\leq\norm{f}_r.
\end{equation}

\medskip
\noindent\emph{Base change of the contraction.}
After base change from $A_0$ to $A_r^\circ$, the block elimination and
contraction of Section~2 remain valid, with the same
contraction constant $\abs{\xi}$.  In particular, the image of $P$ is the
unique fixed point in the corresponding space
$\mathfrak I_k(A_r^\circ)$.  Indeed, units remain units, and the
weighted-column form is preserved under base
change.  The recursive elimination formulas in
Lemma~2.4 therefore make sense over $A_r^\circ$.  In the proof of
Theorem~2.5, the only norm estimate involving $\sigma$
is
\[
  \norm{I^\sigma-I'{}^\sigma}
  \leq\norm{I-I'},
\]
which follows from \eqref{eq:sigma-nonexpanding-small-disc}; equality is not
required.  All remaining estimates are unchanged.  Since
$\mathfrak I_k(A_r^\circ)$ is complete, the contraction has a unique fixed
point, and the image of $P$ is visibly one such fixed point.

The differential equation associated with the connection is
\begin{equation}\label{eq:horizontal-solution-equation}
  \partial(S)=-NS,
  \qquad S(0)=1.
\end{equation}
The standard local existence theorem for analytic differential equations
over a characteristic-zero non-Archimedean field gives, after decreasing
$r$, a unique solution $S\in\operatorname{GL}_n(A_r)$.  Decreasing $r$
once more if necessary, we may and do assume
\[
  \norm{S-1}_r<1,
  \qquad S,S^{-1}\in\operatorname{GL}_n(A_r^\circ).
\]
The horizontal subspace of $M\otimes_A A_r$ is then an
$n$-dimensional $K$-vector space, with basis given by the entries of the
row $\bm eS$ (equivalently, by the vectors represented by the columns of
$S$ in the basis $\bm e$).

\begin{prop}[Block-LU realization of the contraction fixed point]
\label{prop:differential-lu-realization}
Set
\[
  C=S^{-1}BS^\sigma.
\]
Then $C=B(0)\in\operatorname{GL}_n(K)$.  If
$P_0=P(0)$, $Q_0=Q(0)$, and
\[
  \widetilde S=SP_0,
\]
then, after possibly decreasing $r$, $\widetilde S$ has a unique normalized
block LU factorization
\begin{equation}\label{eq:differential-block-lu}
  \widetilde S=JR,
  \qquad
  J\in\mathfrak I_k(A_r^\circ),
\end{equation}
with $R\in\operatorname{GL}_n(A_r^\circ)$ upper block triangular.  Its
lower factor is the restriction of the contraction fixed point over $A_0$:
\[
  J=P|_{A_r}.
\]
Moreover,
\begin{equation}\label{eq:differential-block-conjugacy}
  QR^\sigma=RQ_0.
\end{equation}
\end{prop}

\begin{proof}
From \eqref{eq:horizontal-solution-equation},
\eqref{eq:derivation-frobenius-relation}, and
\eqref{eq:matrix-connection-compatibility}, direct differentiation gives
\[
 \begin{aligned}
  \partial(S^{-1}BS^\sigma)
   &=S^{-1}\bigl(NB+\partial(B)-qt^{q-1}BN^\sigma\bigr)S^\sigma\\
   &=0.
 \end{aligned}
\]
The constants of $A_r$ for $d/dt$ are exactly $K$, so $C\in M_n(K)$.
Because $S(0)=1$, \eqref{eq:differential-frobenius-lift} gives
$C=B(0)$; it is invertible because $B$ is invertible over $A$.

Evaluation of \eqref{eq:differential-fixed-triangularization} at $t=0$
gives
\[
  B(0)P_0=P_0Q_0.
\]
Since $\sigma$ fixes $K$, it follows that
\begin{equation}\label{eq:differential-normalized-solution}
  B\widetilde S^\sigma
  =BS^\sigma P_0
  =SB(0)P_0
  =SP_0Q_0
  =\widetilde S Q_0.
\end{equation}
Now $\widetilde S(0)=P_0$ is lower block unitriangular.  Hence every
leading principal minor associated with a block truncation of
$\widetilde S$ takes the value $1$ at
$t=0$.  After decreasing $r$, all these minors are units of $A_r^\circ$,
so normalized block Gaussian elimination gives the unique factorization
\eqref{eq:differential-block-lu}.

Substitution of $\widetilde S=JR$ into
\eqref{eq:differential-normalized-solution} yields
\[
  BJ^\sigma
  =J\,RQ_0(R^\sigma)^{-1}.
\]
Put $U=RQ_0(R^\sigma)^{-1}$.  This matrix is upper block triangular.
Moreover, because $J\in\mathfrak I_k(A_r^\circ)$, the weighted-column
argument and leading-minor calculation of Lemma~2.3 show
that the $j$th block column of $U=J^{-1}BJ^\sigma$ is divisible by $\xi^j$
and that its normalized diagonal pivots are units.  Thus $J U$ has exactly
the normalized form required in Lemma~2.4.  Uniqueness of that
factorization shows that $J$ is a fixed point of the base-changed
contraction.  The uniqueness established in the base-change argument above
therefore gives $J=P$.
Finally, compare $B\widetilde S^\sigma=\widetilde S Q_0$ with
$BP^\sigma=PQ$ and use $\widetilde S=PR$; cancellation of $P$ gives
\eqref{eq:differential-block-conjugacy}.
\end{proof}

For $0\leq i\leq k$, write $R_i=R_{ii}$ and recall that
$D_i(0)=(Q_0)_{ii}$.  For a block matrix $V$ with invertible preceding
leading block, its $i$th Schur complement is
\[
  \SC_i(V)=V_{ii}-V_{i,<i}V_{<i,<i}^{-1}V_{<i,i},
\]
with $\SC_0(V)=V_{00}$.  We also write
\[
  \widehat\Delta_i(V)=\frac{\Delta_i(V)}{\Delta_{i-1}(V)}.
\]
Left multiplication by a lower block unitriangular matrix does not change a
Schur pivot.  This observation identifies the diagonal blocks of $R$
explicitly in the next theorem.

The lower factor recovers the slope filtration and also relates it to the
connection.  The basis $\bm e\widetilde S$ consists of horizontal elements,
while
\[
  \bm e\widetilde S=\bm eP R.
\]
Because $R$ is upper block triangular with invertible diagonal blocks,
the successive spans of the first block columns agree on both sides.
This also proves that the slope filtration is preserved by the connection
over $A$.  Indeed, its connection matrix in the basis $\bm eP$ is
\[
  N_P=P^{-1}NP+P^{-1}\partial(P)\in M_n(A).
\]
On the smaller disc, the equation for $\widetilde S=PR$ gives
$N_P=-\partial(R)R^{-1}$, which is upper block triangular.  Since
$A\to A_r$ is injective, $N_P$ is upper block triangular over $A$ itself.

\begin{thm}[Exact Schur formula for the graded Frobenius blocks]
\label{thm:differential-schur-formula}
For every $0\leq i\leq k$,
\[
  R_i=\SC_i(\widetilde S),
  \qquad
  \det R_i=\widehat\Delta_i(\widetilde S),
\]
and
\begin{equation}\label{eq:differential-exact-schur}
  D_i
  =\SC_i(\widetilde S)D_i(0)
   \bigl(\SC_i(\widetilde S)^\sigma\bigr)^{-1}.
\end{equation}
The complete expression on the right, formed in $M_{h_i}(K[[t]])$,
belongs to $\xi^i\operatorname{GL}_{h_i}(A_0)$ and equals the globally
defined graded Frobenius matrix $D_i$.
\end{thm}

\begin{proof}
The factorization $\widetilde S=PR$ and invariance of Schur pivots under
left multiplication by a lower block unitriangular matrix give
$R_i=\SC_i(\widetilde S)$ and the determinant identity.  Taking the $i$th
diagonal block in \eqref{eq:differential-block-conjugacy} gives
\[
  D_iR_i^\sigma=R_iD_i(0),
\]
which proves \eqref{eq:differential-exact-schur} on the disc of
convergence of $S$.  Both sides have the same Taylor series at $t=0$.
Since $D_i\in\xi^i\operatorname{GL}_{h_i}(A_0)$, the assertion about the
complete formal expression follows.
\end{proof}

Equivalently, the last factor in \eqref{eq:differential-exact-schur} may be
written as $\SC_i(\widetilde S^\sigma)^{-1}$.  Thus the full Schur
complement determines the diagonal block of the upper factor, whereas the
leading principal minors determine its determinant.  The normalization is
essential: the matrix being factorized is $\widetilde S=SP_0$.  If $B(0)$
is already upper block triangular, then $P_0=1$ and $\widetilde S=S$.

The matrices $S$ and $R_i$ need not extend individually to the closed
unit disc.  Expressions involving them will first be formed in $K[[t]]$.
Once a complete expression is shown to belong to $A_0$, it can be
specialized at any bounded $\mathcal O_K$-algebra homomorphism from $A_0$.
Such a specialization refers to the complete expression, without requiring
separate evaluations of its local factors.

\begin{cor}[Periodic fibre characteristic factors]
\label{cor:differential-periodic-factors}
For $d\geq1$, form the matrix
\[
  G_{i,d}
  =\SC_i(\widetilde S)D_i(0)^d
   \bigl(\SC_i(\widetilde S)^{\sigma^d}\bigr)^{-1}
\]
in $M_{h_i}(K[[t]])$.  Then
\begin{equation}\label{eq:differential-periodic-factor}
  G_{i,d}
  =D_iD_i^\sigma\cdots D_i^{\sigma^{d-1}}
  \in\xi^{id}\operatorname{GL}_{h_i}(A_0).
\end{equation}
Let $K_x/K$ be a finite valued extension and let
\[
  x:A_0\longrightarrow\mathcal O_{K_x}
\]
be a bounded $\mathcal O_K$-algebra homomorphism satisfying
$x\circ\sigma^d=x$.  The local Frobenius matrix on the $i$th graded
fibre and its characteristic polynomial are
\[
  \Phi_{i,x}=G_{i,d}(x),
  \qquad F_{i,x}(T)=\det\bigl(1-TG_{i,d}(x)\bigr).
\]
With the valuation normalized by $v_x(\xi)=1$, this is the complete
slope-$i$ factor of the fibre characteristic polynomial when slopes are
measured using $v_x/d$.
\end{cor}

\begin{proof}
Apply $\sigma^s$ to \eqref{eq:differential-exact-schur} and multiply in
order for $0\leq s<d$.  Since $\sigma$ fixes $K$, adjacent factors cancel
to give \eqref{eq:differential-periodic-factor} in $M_{h_i}(K[[t]])$.
The product of the global matrices $D_i^{\sigma^s}$ proves the asserted
membership in $\xi^{id}\operatorname{GL}_{h_i}(A_0)$.  Specialization
gives the fibre formula.  Finally, $P\in\operatorname{GL}_n(A_0)$ and
periodicity identify the local Frobenius with an upper block triangular
matrix whose diagonal blocks are the $G_{j,d}(x)$.  The normalized block
$\xi^{-jd}G_{j,d}(x)$ is integral and has integral inverse, so all its
eigenvalues have valuation zero.  Thus the $j$th block has slope $j$,
which proves the slope-factor assertion.
\end{proof}

\begin{cor}[Rank-one leading-minor and root formulas]
\label{cor:differential-rank-one}
Suppose $h_i=1$.  For $d\geq1$, define the formal power series
\[
  u_{i,d}
  =\frac{\widehat\Delta_i(\widetilde S)}
         {\widehat\Delta_i(\widetilde S)^{\sigma^d}}
  =\frac{
     \Delta_i(\widetilde S)\Delta_{i-1}(\widetilde S)^{\sigma^d}
   }{
     \Delta_{i-1}(\widetilde S)\Delta_i(\widetilde S)^{\sigma^d}
   }.
\]
Then $u_{i,d}\in A_0^\times$, and
\begin{equation}\label{eq:differential-rank-one-schur}
  D_i=D_i(0)u_{i,1}
  =D_i(0)
   \frac{\widehat\Delta_i(\widetilde S)}
        {\widehat\Delta_i(\widetilde S)^\sigma}.
\end{equation}
For a periodic specialization $x:A_0\to\mathcal O_{K_x}$ as in
\Cref{cor:differential-periodic-factors}, the Frobenius eigenvalue on
this rank-one graded fibre is
\begin{equation}\label{eq:differential-rank-one-root}
  \gamma_{i,x}
  =\prod_{s=0}^{d-1}D_i^{\sigma^s}(x)
  =D_i(0)^d u_{i,d}(x).
\end{equation}
\end{cor}

\begin{proof}
When $h_i=1$, the identity
$\SC_i(\widetilde S)=\widehat\Delta_i(\widetilde S)$ turns
\eqref{eq:differential-exact-schur} into
\eqref{eq:differential-rank-one-schur}.  By
\eqref{eq:differential-periodic-factor},
\[
  u_{i,d}
  =\prod_{s=0}^{d-1}
       \left(\frac{D_i}{D_i(0)}\right)^{\sigma^s}.
\]
Each factor belongs to $A_0^\times$, since
$D_i\in\xi^iA_0^\times$ and
$D_i(0)\in\xi^i\mathcal O_K^\times$.  Therefore $u_{i,d}\in A_0^\times$,
and specialization gives \eqref{eq:differential-rank-one-root}.
\end{proof}

In \eqref{eq:differential-rank-one-root}, $u_{i,d}(x)$ is obtained by
specializing the complete leading-minor quotient, whose Taylor series
defines an element of $A_0^\times$.  Periodicity applies to elements of
$A_0$.  It cannot be used to cancel the local factors
$\widehat\Delta_i(\widetilde S)$ and
$\widehat\Delta_i(\widetilde S)^{\sigma^d}$ at $x$, since they need not
be individually defined there.

Taking determinants in \eqref{eq:differential-exact-schur} gives
\[
  \det D_i
  =\det D_i(0)
   \frac{\widehat\Delta_i(\widetilde S)}
        {\widehat\Delta_i(\widetilde S)^\sigma}.
\]
In higher rank this records only the determinant of $D_i$; the full block in
\eqref{eq:differential-exact-schur} requires the entire Schur complement,
equivalently its bordered-minor entries.

Together with \Cref{cor:differential-periodic-factors}, the global
identity \eqref{eq:differential-periodic-factor} determines the complete
characteristic polynomial, and hence the multiset of Frobenius eigenvalues,
on each graded fibre.  In higher rank it does not canonically label those
eigenvalues; in rank one, \eqref{eq:differential-rank-one-root} gives the
individual eigenvalue.  Reciprocal zeros and poles of an $L$-function are
then obtained through an applicable trace formula, with its slope shifts
and possible cancellations.

The differential equation thus gives exact Schur expressions for the
graded Frobenius matrices defined over $A_0$.  Although the solution
matrix is constructed locally, the complete expressions in
\Cref{thm:differential-schur-formula,cor:differential-periodic-factors}
belong to $A_0$ and can be specialized at its periodic points.
\Cref{cor:differential-rank-one} expresses the corresponding rank-one
eigenvalues through leading principal minors.  These exact expressions
complement the asymptotic formulas obtained from the contraction iterates
in Section~2.

This example illustrates how additional differential structure may be used
when a finite-rank $\sigma$-module arises at the cohomological level.
The chain-level application remains the main subject of this paper.
Investigating such differential structures in broader cohomological or
multi-parameter settings lies beyond its scope.

\section*{Acknowledgements}
The author thanks Prof. Steven Sperber  and Prof. Daqing Wan for communication on an earlier version of the manuscript. Their suggestions greatly improved the presentation of the paper. The author acknowledges the assistance of ChatGPT  during the development and preparation of this paper. In the course of discussions, ChatGPT suggested expressing the finite Frobenius approximations in terms of Schur complements. The author verified this observation, developed it into the asymptotic formulas proved here, and incorporated it into the mathematical framework. ChatGPT also assisted in organizing the arguments, refining the exposition, standardizing notation, checking references and cross-references, and preparing the \LaTeX{} manuscript.

The research problem, the contraction construction, its applications to slope factors of $L$-functions, and the overall chain-level and cohomological framework were conceived and developed by the author. The author takes full responsibility for the entire content of the paper.


\begin{thebibliography}{99}

\bibitem{AS89}
A.~Adolphson and S.~Sperber,
\emph{Exponential sums and Newton polyhedra: cohomology and estimates},
Ann. of Math. (2) \textbf{130} (1989), 367--406.

\bibitem{AS12}
A.~Adolphson and S.~Sperber,
\emph{On unit root formulas for toric exponential sums},
Algebr. Number Theory \textbf{6} (2012), 573--585. 

\bibitem{ASZ} 
A.~Adolphson, S.~Sperber and H.~Zhang, 
\emph{Maximal unipotent monodromy and $p$-adic eigenvalue formulas}, preprint, arxiv: https://arxiv.org/pdf/2611.05886.

\bibitem{BGR}
S.~Bosch, U.~G\"untzer, and R.~Remmert,
\emph{Non-Archimedean Analysis: A Systematic Approach to Rigid Analytic
Geometry}, Grundlehren Math. Wiss. 261,
Springer-Verlag, Berlin, 1984.

\bibitem{De80}
P.~Deligne,
\emph{La conjecture de Weil. II},
Publ. Math. IH\'ES \textbf{52} (1980), 137--252.

\bibitem{DL91}
J.~Denef and F.~Loeser,
\emph{Weights of exponential sums, intersection cohomology, and Newton
polyhedra}, Invent. Math. \textbf{106} (1991), 275--294.

\bibitem{Dw60}
B.~Dwork,
\emph{On the rationality of the zeta function of an algebraic variety},
Amer. J. Math. \textbf{82} (1960), 631--648.

\bibitem{Dw69}
B.~Dwork,
\emph{$p$-adic cycles},
Publ. Math. IH\'ES \textbf{37} (1969), 27--115.

\bibitem{Dw73}
B.~Dwork,
\emph{On $p$-adic differential equations. IV: generalized hypergeometric
functions as $p$-adic analytic functions in one variable},
Ann. Sci. \`Ecole Norm. Sup. (4) \textbf{6} (1973), 295--316.

\bibitem{Dw74}
B.~Dwork,
\emph{Bessel functions as $p$-adic functions of the argument},
Duke Math. J. \textbf{41} (1974), 711--738.

\bibitem{HS14}
C.~Haessig and S.~Sperber,
\emph{L-functions associated with families of toric exponential sums},
J. Number Theory \textbf{144} (2014), 422--473.

\bibitem{Sp80}
S.~Sperber,
\emph{Congruence properties of the hyperkloosterman sum},
Compositio Math. \textbf{40} (1980), 3--33.

\bibitem{W93}
D.~Wan,
\emph{Newton polygons of zeta functions and $L$-functions},
Ann. of Math. (2) \textbf{137} (1993), 249--293.

\bibitem{W99}
D.~Wan,
\emph{Dwork's conjecture on unit root zeta functions},
Ann. of Math. (2) \textbf{150} (1999), 867--927.

\bibitem{W00a}
D.~Wan,
\emph{Higher rank case of Dwork's conjecture},
J. Amer. Math. Soc. \textbf{13} (2000), 807--852.

\bibitem{W00b}
D.~Wan,
\emph{Rank one case of Dwork's conjecture},
J. Amer. Math. Soc. \textbf{13} (2000), 853--908.

\bibitem{W04}
D.~Wan,
\emph{Variation of $p$-adic Newton polygons for $L$-functions of exponential
sums}, Asian J. Math. \textbf{8} (2004), 427--472.

\bibitem{YZ26}
L.~Yang and H.~Zhang,
\emph{Unit roots of the unit root $L$-functions},
Math. Ann. \textbf{394} (2026), article 28.

\bibitem{Zh}
H.~Zhang,
\emph{Formulas for roots of $L$-functions of the Airy exponential sums},
Finite Fields Appl. \textbf{49} (2018), 29--48.

\end{thebibliography}
\end{document}